\documentclass[letterpaper]{article} 
\usepackage{aaai24}  
\usepackage{times}  
\usepackage{helvet}  
\usepackage{courier}  
\usepackage[hyphens]{url}  
\usepackage{graphicx} 
\usepackage{natbib}  
\usepackage{caption} 
\usepackage{algpseudocode}
\usepackage{algorithm}
\usepackage{amsmath}
\usepackage{amssymb}
\usepackage{amsthm}
\usepackage{amsfonts}
\usepackage{adjustbox}
\usepackage{enumerate} 
\usepackage[shortlabels]{enumitem}
\usepackage{mathtools}
\usepackage{cite}
\usepackage{xcolor}
\usepackage{multicol}

\usepackage{subcaption}
\usepackage{multicol}
\usepackage{graphicx}

\algnewcommand{\IIf}[1]{\State\algorithmicif\ #1\ \algorithmicthen}
\algnewcommand{\EndIIf}{\unskip\ \algorithmicend\ \algorithmicif}
\def\A{\mathcal{A}}
\def\X{\mathcal{X}}
\def\Y{\mathcal{Y}}

\def\D{\mathcal{D}}
\def\N{\mathcal{N}}
\def\M{\mathcal{M}}
\def\L{\mathcal{L}}

\def\B{\mathcal{B}}

\def\P{\mathcal{P}}

\def\I{\mathcal{I}}
\def\J{\mathcal{J}}
\def\T{\mathcal{T}}
\def\O{\mathcal{O}}
\def\Z{\mathcal{Z}}

\def\calP{\mathcal P}

\def\calN{\mathcal N}

\def\bd{{\mathbf d}}

\def \bp {{\mathbf{p}}}
\def \bq {{\mathbf{q}}}

\def\argmax{\mathop{\rm arg\,max}\limits}

\newcommand{\commentbyKG}[1]{}

\newcommand{\commentbyET}[1]{}

\providecommand{\keywords}[1]{\textbf{\textit{Index terms---}} #1}

\newcommand{\subsubsubsection}[1]{\paragraph{#1}\mbox{}}
\newtheorem{theorem}{Theorem}[section]
\newtheorem{lemma}[theorem]{Lemma}

\newtheorem{definition}[theorem]{Definition}
\numberwithin{equation}{section}

\title{INVALS: A Forward Looking Inventory Allocation System}
\author {
    Shiv Krishna Jaiswal, 
    Karthik S. Gurumoorthy, 
    Etika Agarwal, 
    Shantala Manchenahally
}
\affiliations {
    Walmart Global Tech, Bangalore, India\\
    {\{shivkrishna.jaiswal, karthik.gurumoorthy, etika.agarwal, shantala.manchenahally\}@walmart.com}
  
}

\begin{document}

\maketitle
\begin{abstract}
    We design an \textbf{Inv}entory \textbf{Al}location \textbf{S}ystem (INVALS) that, for each item-store combination, plans the quantity to be allocated from a warehouse that replenishes multiple stores using trailers, while respecting typical operational constraints. We formulate a linear objective function which, when maximized, determines the allocation plan by considering not only the immediate store needs, but also its future (forward) expected demand. This forward-looking allocation significantly improves the utilization of labor and trailers in the warehouse. To reduce overstocking, we adapt from our objective to prioritize allocating those items in excess which are sold faster at the stores, keeping the days of supply (DOS) to a minimum. For the proposed formulation, which is an instance of Mixed Integer Linear Programming (MILP), we present a scalable algorithm using the concepts of submodularity and optimal transport theory by: (i) sequentially adding trailers to stores based on maximum incremental gain, (ii) transforming the resultant linear program (LP) instance to an instance of capacity constrained optimal transport (COT), solvable using double entropic regularization and incurring the same computational complexity as the Sinkhorn algorithm. Compared against the planning engine that only allocates for immediate store needs, INVALS increases labor utilization by $34.70\%$ and item occupancy in trailers by $37.08\%$ on average. The DOS distribution is also skewed to the left, indicating that higher-demand items are allocated in excess, reducing the days they are stocked. We empirically observed that for $\approx$ 90\% of the replenishment cycles, the allocation results of INVALS are identical to the globally optimal MILP solution.
\end{abstract}


\section{Introduction}
A retail giant typically consists of thousands of stores spread over a vast geographic area, with each store offering tens of thousands of products. Each store receives its replenishment regularly from the warehouse to which it is assigned. For every combination of (item, store) the forecasting models predict the expected demand for multiple days into the future. Using these forecast values and the current inventory on-hand, each store for every item determines the quantity it would require to meet the customer demand for those days in the future. This is the store demand a.k.a. need sent to the warehouse which, based on the requirements received from multiple stores, generates a plan to allocate the available inventory--- the \emph{inventory allocation plan}. Once an allocation plan of items to store is created, shipments from the warehouse are moved to stores with the help of trailers. Trailers are dispatched to stores only on certain days of the week based on its replenishment cycle and a store may not receive trailers every day. The time between two consecutive receiving of trailers is called the \emph{coverage period}. The granularity of allocation is a box level known as \textit{warehouse-packs} (whpacks) consisting of a fixed number of same item. A stock-out at a store occurs if the inventory for a particular product at the store goes to zero. This is undesirable as it not only leads to lost sales but also poor customer experience. Holding costs are incurred by stores to maintain on-hand inventory of products. These are comprised of quantities such as electricity costs for the store, refrigeration costs for food items, and storage area maintenance costs etc., and increase with the amount of inventory kept on-hand. As the warehouse has limited inventory and labour availability, allocation of items must be carefully done and prioritised to stores requiring them the most.

 The store needs comprises of two parts: (i) Coverage/immediate need: Quantity required for current replenishment cycle that needs to be shipped to prevent stock-outs at the stores, and (ii) Pull-forward need: Need on the next few days in the future after the coverage period. While it is not necessary to allocate items more than the immediate need for the current replenishment cycle, doing so can prove beneficial on multiple aspects. As an example, it is often seen that store demand for Mondays is generally low and peaks during weekends. If the inventory planned for Monday only caters to the Monday's need, then the workforce available at the warehouse will be under-utilized during this planning cycle. As fewer items are loaded to the trailer, its capacity utilization will be equivalently low resulting in higher per-item transportation cost. A contrasting picture emerges during the weekend planning where the warehouse lacks manpower to process the higher demand from the stores, leading to even the coverage need being unmet and resulting in stock-outs. 
 
 In this work, we argue that the planning system should be able to consider requirements from upcoming future days while deciding the replenishment for the replenishment cycle, provided no business constraints are violated. We term this additional allocation from future demand as \emph{Pull Forward} (PF). As we show in our experiments, PF leads to higher labour usage and improved trailer utilization. 

\subsection{Contributions}
Following is a brief summary of our contributions.
\begin{enumerate}
\item Design of a new forward looking inventory allocation system (INVALS), that decides to allocate items meeting future demand in addition to fulfilling the immediate store needs, for better and efficient utilization of labour and trailer at the warehouse.
\item Representing INVALS from the lens of an optimization framework which trades off between allocating more items to stores, and simultaneously minimising the storage of excessive inventory for longer duration of time.
\item Development of a fast algorithm using the concepts of submodularity and optimal transport theory to solve the optimization framework. We demonstrate that the resultant linear program (LP) instance can be transformed to an instance of capacity constrained optimal transportation (COT) by inclusion of multiple pseudo variables. This transformation process could be of independent interest and open a separate thread of research in identifying the structure of LP instances amenable to such reductions. This assumes significance given that LP formulations appear in plethora of problems studied in operations research.
\end{enumerate}

\section{Replenishment Systems and Related Work}
The replenishment algorithms generally studied in the literature identify the best replenishment strategy for a single item, to maintain proper inventory levels under different settings~\cite{Bosman2022}. Some of the popular methods include: (i) 
\emph{Reorder point} where the current inventory level triggers the action to replenish the item, (ii) \emph{Periodic} where the inventory is replenished at specific time intervals at the predetermined review points, (iii) \emph{Actual Demand} where the item is restocked as and when the demand arises, and (iv) \emph{Top off} where restocking is done to fixed levels during business lean time~\cite{Khan99}. Though our work is based on periodic replenishment strategy, determining the allocation plan simultaneously for 10000's of products across multiple stores satisfying a plethora of constraints is a novel problem, and to the best of our knowledge has received very little attention. Approaches in~\cite{Martino2016} and the references there in are specific to fashion retail industry, where the replenishment problem is to maximize profit ---defined as difference between revenue and incurred costs--- satisfying a whole market demand that follows a uniform distribution, and subject to a budget constraint. Replenishment of multiple items is studied as a Joint Replenishment Problem (JRP)~\cite{Muriel2022, Goyal1989, Khouja2008, Peng2022, Zhuang2023} when certain fixed costs are incurred regardless of mix or quantities of the items ordered. The objective of these approaches is usually to minimize a total monetary cost while satisfying demand. In contrast, our work does not have an explicit monetary cost minimisation objective. We introduce a new concept of pulling-forward items from future demand into the current replenishment cycle to maximize labour and trailer usage.

\subsection{Existing Replenishment System}
The planning system deployed in large retail systems are rule-based in nature, drafting an allocation plan that only focuses on the coverage need, and considers only the inventory availability at the warehouse. With lack of visibility to labour or trailer usage, these allocations plans often turn out to be infeasible to execute by ground operations. The reasons could be shortage of labour workforce, exceeding the trailer capacity to transport the allocated items, and in many cases severe under utilization of trailer capacity. To counter this, multiple downstream systems update this plan, each making it feasible only with respect to a specific criteria. Operating sequentially based on rules, these systems lack visibility to the overall replenishment system and introduce their own quantum of sub-optimality to labour and trailer usage. Further, maintaining multiple such systems is inefficient and costlier.

\subsection{INVALS}
\label{sec:INVALS}
INVALS is a forward-looking planning engine that creates the inventory allocation plan by not only looking into the immediate store needs, but also multiple days into the future while optimizing for labour and trailer usage, reducing out of stock, and minimizing excess inventory. By maximizing the tailor-made linear objective function \eqref{eq:objective} stated in Sec.~\ref{sec:costfunction} and incorporating various business constraints described in Sec.~\ref{sec:constraints}, INVALS solves a unified problem leading to improved operational efficiencies and making available the right items needed at the right time at the right store.

The primary challenge in designing INVALS is at the formulation of the right objective function whose maximisation leads to higher labour and trailer usage by pulling-forward items, and also ensuring that the coverage need for any store is always prioritised before allocating items for future need of other stores to reduce stock-outs. Allocating more items can have the adverse effect of increased inventory at the stores and ensuing holding costs. To reduce these costs, our formulation should allocate only those surplus items to stores which will be sold faster, keeping the number of days the items last, known as days of supply (DOS), at the store to a minimum. The framework should have the flexibility to accommodate business constraints where trailers can be prioritised to stores with longer replenishment cycles, and guarantee that trailers are loaded with a minimum quantity of items to justify their transportation cost. The second challenge is in developing a fast optimization algorithm respecting the standard supply chain constraints, the solution to which determines the inventory allocation plan. 
\section{Mathematical Formulation}
\label{sec:formulation}
\subsection{Notations}
\label{sec:notations}
Let $\I$ and $\J$ denote the set of items and stores respectively indexed by $i$ and $j$. Often items are categorized into multiple channels such as those that can and cannot be moved via a conveyor belt etc. We denote by $\L$ the set of such available channels. For $l \in \L$,  $\I_l$ represents items belonging to channel $l$, and $H_l$ denotes the available labour capacity in terms of number of whpacks that can processed per day at the warehouse for items in $\I_l$. We represent the trailer max-capacity and min-capacity by $M$ and $m$, the shelf-capacity of the $i^{th}$ item in store $j$ by $C_{ij}$, and the maximum number of trailers that can be dispatched to store $j$ by the integer $R_j$. Coverage period is denoted by time $t = 0$, and time periods for future need of the $i^{th}$ item on $j^{th}$ store by $t_{ij}$. Based on the store needs  $D_{ij}^{t}$ for each time $t$ in $\T_{ij} = \{0\} \cup \{1,2, \ldots t_{ij}\}$ for $i \in \I, j\in \J$, INVALS does the inventory allocation to satisfy the immediate need at $t=0$ and perform PF by looking forward into the needs of the future $t_{ij}$ days. Let $S_i$ denote the available inventory of the $i^{th}$ item at the warehouse, $p_j$ represent the store-priority, and $q_{ij}^t$ indicate the item-store priority for sending $i^{th}$ item to $j^{th}$ store on $t^{th}$ day. The quantities $p_j$ and $q_{ij}^t$ are discussed below in Sec.~\ref{sec:costfunction}. The hyper-parameters are denoted by $\alpha [t], t \in \T_{ij}$ which is an exponentially decreasing function of $t$, and large positive constants $\beta$ and $\gamma$ with $\gamma > \beta$. 

Our aim is to determine the variables: (i) $\mathcal{D} = \{d_{ij}^t\}$ indicating the quantity of $i^{th}$ item allocated to $j^{th}$ store corresponding to the $t^{th}$ day, and (ii) $\X = \{x_j\}$ representing the number of trailers dispatched to each store $j$. For simplicity, we introduce additional variables $s_{ij} = \sum_{t \in \T_{ij}} d_{ij}^t$ denoting the allocation of item $i$ to store $j$ over all time periods, $s_j = \sum_{i \in \I}s_{ij}$ representing the total inventory allocation to store $j$, $y_j = min(1, x_j)$ which indicates whether any trailer is assigned to store $j$, and $b_j$ corresponding to the extent of trailer min-capacity violation in the last trailer dispatched to the $j^{th}$ store as explained below. We define the set $\B = \{b_j\}$ comprising of the min-capacity breach slack variables.

\subsection{Utility Function}
\label{sec:costfunction}
We propose to \emph{maximize} the following utility function:
\begin{align}
\label{eq:objective}
    Obj(\D, \X, \B) &= \underbrace{\left(\sum_{\substack{i \in \I,j\in \J\\t\in \T_{ij}}}\alpha[t]q_{ij}^t d_{ij}^{t}\right)}_{\text{Item-store allocation utility}} + \underbrace{\beta \left(\sum_{\substack{j \in \J}} p_{j} x_{j}\right)}_{\text{Store prioritization utility}} \nonumber \\
    &\underbrace{- \gamma \left(\sum_{\substack{j \in \J} }b_j \right)}_{\text{Trailer min-capacity cost}}
\end{align}

subject to the various constraints described in Sec.~\ref{sec:constraints}. Our objective comprises of three terms namely: (i) Item to store allocation utility, (ii) Store prioritization utility, and (iii) Cost due to LTMC (less than min capacity) trailers, each described and motivated below.

The first term is the \emph{allocation utility} made-up of three quantities: (i) Allocation ($d_{ij}^t$) of $i^{th}$ item to $j^{th}$ store for $t^{th}$ day, (ii) Quantification of importance of the $i^{th}$ item to $j^{th}$ store at time $t$ ($q_{ij}^t$) based on how fast the item is sold at the store as explained in Sec.~\ref{sec:experiments}, and (iii) Exponentially decreasing function $\alpha[t]$ to prioritise allocating the immediate need for any store over the future need of others. Sans this function, the plan could result in surplus allocation of an item (more than the coverage period need) at one store and stock-out (quantity less than the coverage-period need) of the same item on another store.

The middle term is the \emph{store prioritization utility}. As trailers are dispatched to stores only on specific days of the week, not sending items to stores with longer coverage period has higher negative impact on sales compared to those which receive trailers more frequently. It is also important from the business perspective that stores opened in the recent past are given more priority. To model the relative ranking among stores, we introduce the store prioritization term $p_j \in \mathbb{N}$ representing how critical is it for the store to receive a trailer as INVALS is planning the item allocation. The higher the value of $p_j$, greater the importance of assigning trailers to the store. We control the impact of this term with a parameter $\beta > 0$.

The last term refers to the \emph{Less Than Min Capacity (LTMC)} trailer cost appearing in the negative. INVALS encourages to utilize the trailer capacity maximally. However, due to scarcity of labour and inventory at warehouse, or insufficient need from store, trailers may not be filled completely resulting in higher transportation cost per item. Hence we desire a minimum shipment quantity $m$ that should be loaded on each trailer assigned to a store. Nevertheless, at times it may still be required for the trailer to be dispatched with less than $m$ shipments for smooth customer buying experience at the store. To avoid infeasibility in the formulation, we introduce a slack variable $b_j$ to represent min-capacity breach of last trailer headed towards $j^{th}$ store. By defining a parameter $\gamma$, we control the extent to which min-capacity breach can be allowed. Setting $\gamma = \infty$ is equivalent to disallowing min-capacity breach, and any other non-negative $\gamma$ implies that the trailer min-capacity can be breached provided the total benefit of allocating inventory using this trailer is at least $\gamma$.

\subsection{Constraints}
\label{sec:constraints}
The following constraints must be respected for the allocation plan to be operationally feasible.

\noindent \textbf{Inventory Constraint}
For any item $i$, the sum of its allocation across all stores should be within the inventory available at warehouse, i.e., $\sum_{j \in \J}s_{ij} \leq S_i$. 

\noindent \textbf{Warehouse Labour Constraint}
For every labour channel $l$ the total processed quantity of items requiring that labour type should be within the available workforce $H_l$, represented by the constraint $\sum_{i \in \I_l,j \in \J} s_{ij} \leq H_l, \forall l \in \L$.

\noindent \textbf{Planned trailer constraint}
The condition $x_j \leq R_j, x_j \in \mathbb{N} \cup \{0\}, \forall j \in \J$ bounds the number of trailers $x_j$ dispatched to store $j$ to be within a pre-specified constant $R_j$. This condition is related to the availability of store labour to unload these items upon arrival.

\noindent \textbf{Trailer max-capacity constraint}
The condition $s_j \leq Mx_j, \forall j \in \J$ ensures that total items shipped to $j^{th}$ store is within the available trailer capacity.

\noindent \textbf{Trailer min-capacity constraint}
As mentioned before, we desire that each trailer is loaded with at least $m$ whpacks for transportation costs to be justified. However, sometimes the allocation plan may breach the min-capacity and dispatch a trailer with items less than $m$ whpacks due to labour or inventory shortage at the warehouse, or inadequate store need. If multiple trailers are dispatched to a store, the min-capacity breach can only occur in the last trailer, as all earlier trailers have to be fully loaded to their max-capacity $M$ before assigning an additional trailer. We represent the same by introducing a non-negative slack variable $b_j \leq m$ for each store and enforce that total store allocation $s_j$ should satisfy $s_j \geq M(x_j-y_j) + my_j -b_j$. The desirable solution is to have $b_j=0, \forall j$ and avoid min-capacity breaches.

\noindent \textbf{Shelf capacity constraint}
We impose $s_{ij} \leq C_{ij}, \forall i\in \I, j\in \J$ to make sure that the total (item, store) allocation $s_{ij}$ does not exceed the available shelf-capacity $C_{ij}$ at the store to hold the inventory upon arrival. 

\noindent \textbf{Demand constraint}
The condition $0 \leq d_{ij}^{t} \leq D_{ij}^{t}, \forall i\in \I, j\in \J,t\in \T_{ij}$ ensures that allocation for any day $t$ does not exceed the store need $D_{ij}^t$. Ideally, one should enforce that the allocations $d_{ij}^t$ and the min-capacity breach slack variables $b_j$ also satisfy the integrality constraint and assume only non-negative integer values. This would not scale as the number of such decision variables are typically in millions for a large supply chain network. We choose to relax this constraint as our problem structure inherits some properties of total unimodularity (TU)~\cite{Hoffman2010}. As stated below in \eqref{problem:linear_program}, the item allocation problem involving $d_{ij}^t$ is an instance of LP once the variables $x_j$ are fixed, and can be succinctly expressed as $\max_{\bd}\{\bp^T \bd | A\bd \leq \bq\}$, for vectors $\bd, \bp ,\bq$ and matrix $A$. The entries of $A$ are either $\{0,1,-1\}$ and $\bq$ is integral. This does not necessarily guarantee that $A$ is a TU matrix to endow the LP instance with integral optima, where the integrality constraints on the variables $d_{ij}^t$ and $b_j$ are implicitly met. However, in our experiments we observed that for an overwhelming $99.99\%$ of the variables $d_{ij}^t$, the solution is integer valued even without the integrality constraint. Theoretical analysis of this relaxation is beyond the scope of this work.

\section{Proposed Solution}
The formulation in~\eqref{eq:objective} is an instance of MILP and can be solved using any off-the-shelf solver. The presence of integer variables $x_j$ makes the problem NP-Hard, and obtaining the global optimum as the solution to MILP instance is not run-time efficient and increases exponentially in the number of such variables. However, once the number of trailers assigned to each store namely the variables $x_j$ are determined, solving for the $d_{ij}^t$ variables to compute the item allocation to stores reduces to an instance of LP. Let $\X^* = \{x_1^*,x_2^*,\ldots,x_J^*\}$ be the set of store-trailer mappings with entries corresponding to the number of trailers assigned to each store respecting the condition that $x_j^* \leq R_j$. Recall that $y_j^* = \min(1,x_j^*)$. The linear program objective is:
\begin{equation}
\label{problem:linear_program}
 \mathbb{P}_1: g(\X^*) = \max\limits_{d_{ij}^t, b_j | \X^*} \left(\sum\limits_{i,j,t}\alpha[t]q_{ij}^t d_{ij}^{t}\right) - \gamma \left(\sum\limits_{j}b_j \right)
\end{equation}
subject to all the constraints detailed in Sec.~\ref{sec:constraints} barring the planned trailer constraint which depends only on $x_j^*$. We assume $g(\X^*) = 0$ when the feasible set is empty. Bearing this in mind, we develop an incremental iterative algorithm leveraging the concepts of submodularity and optimal transport theory. In this approach we incrementally grow the value $x_j$ of the number of trailers dispatched to any store $j \in \mathcal{J}$ by posing the optimization problem as a selection of best store assignments for a sequence of trailers. Our methodology consists of three steps namely: (i) Determine the store to assign the $n+1^{th}$ trailer given the assignment $\X^n$ of the first $n$ trailers, (ii) Efficiently compute the objective value for the item allocation problem $\mathbb{P}_1$ in ~\eqref{problem:linear_program} for each candidate store $k$ in the previous step, by transforming the LP instance into an instance of capacity constrained optimal transport and using the double regularization method (DRM)~\cite{Wu2022}, and (iii) Once the final number of trailers to be dispatched to each store are determined, solve the item allocation for~\eqref{problem:linear_program} using any standard LP solver just for this one last instance to compute the item allocation $d_{ij}^t$. As DRM is only an approximate algorithm, the usage of LP solver for the last step yields the optimal solution for $d_{ij}^t$ given the final trailer list. \commentbyKG{we empirically observed $99.99\%$ of $d_{ij}^t$ to be integer valued. it may not necessarily yield such overwhelming integral solution.}We describe these steps below. 

\subsection{Incremental Trailer Assignment}
\label{sec:trailerassignment}
The aim of this module is to find the number of trailers assigned to each store. Let $\left(\D^*_{\X}, \B^*_{\X}\right)$ be the optimal point for~\eqref{problem:linear_program}.  The objective in~\eqref{eq:objective} can be rewritten as: 
\begin{equation}
\label{eq:setfunction}
 Obj\left(\D^*_{\X}, \X, \B^*_{\X}\right) = f(\X) = g(\X) + \beta*h(\X)
\end{equation}
where $g(\X)$ is defined in~\eqref{problem:linear_program} and $h(\X)=\sum_{j \in \J}p_jx_j$. Our goal is to maximize $f(\X)$ subject to the planned trailer constraint $x_j \leq R_j, \forall j \in \J$. The function $g(\X)$ satisfies the property of submodularity and non-negativity for chosen values of the hyper-parameters $\alpha[t]$, $\beta$ and $\gamma$, and $h(\X)$ is a non-negative, modular function. It then follows that $f(.)$ is a non-negative, submodular function. These concepts and discussion of submodularity is detailed in Sec.~C of the supplementary.

Submodularity implies diminishing returns where the incremental gain in adding a new element to a set $\mathcal{A}$ is at least as high as adding to its superset $\mathcal{B}$~\cite{fujishige05, lo83}. In our context, it implies that it is more beneficial to assign a trailer to a store when there are fewer trailers assigned. Maximisation of submodular set functions is a well studied problem. As they are $\calN \calP$-hard in the general form, a variety of algorithms have been proposed to find approximate solutions to submodular optimization problems. One of the most popular category of algorithms are the variants of incremental selection of a set using greedy approaches~\cite{Nemhauser78, Buchbinder2014, Buchbinder2017, Kuhnle2019, Sakaue2020}. The methods in~\cite{Buchbinder2014, Kuhnle2019, Sakaue2020} provide approximation guarantees when the submodular function in~\eqref{problem:linear_program} is only non-negative, without necessitating it to be monotone.

A faster variation of the greedy algorithm called the \textit{stochastic greedy algorithm}~\cite{mirzasoleiman15a} progresses as follows. Let $\X^n = \{x_1^n,x_2^n,\ldots,x_J^n\}$ be the number of trailers assigned to each store at the end of iteration $n$. The property of diminishing returns enables us to assign trailers in an incremental fashion, where starting from the set $\X^0=\{0,0,\ldots,0\}$, in every iteration $n+1$, given the allocation of stores for the first $n$ trailers, we identify the best store $k_{n+1}$ for the next trailer from a random subset $\J^{n+1}$ of eligible candidate stores $\A^{n+1} = \{j \in \J: x_j^n < R_j\}$. The selection is made to maximize the incremental gain (IG) $f_{\X^n}(k) = f\left(\X^n \uplus \{k\}\right) - f(\X^n)$. The set $\X^{n}$ is grown incrementally to $\X^{n+1}=\X^n \uplus \{k_{n+1}\}$ by including the chosen store $k_{n+1}$ with the highest IG. We slightly abuse the notation and define $\X^{n+1}=\X^n \uplus \{k\}$ as the set of store-trailer mappings with $x_k^{n+1} = x_k^n+1$, and $x_j^{n+1} = x_j^n, \forall k \not=j$. This way of incrementally selecting the stores eliminates the need of integrality constraints, and is a major factor for reducing the computational complexity. 

\subsubsection{Batching}
For large value of $\beta$, stores in the candidate set $\A^{n+1}$ with higher store-priority value $p_j$ are most likely to yield the highest IG because of the factor $\beta*p_j$.  To this end, we create a batch $\P_p = \{j \in \J: p_j=p\}$ by grouping stores with same store-priority $p$ and in every iteration $n$ try assigning the next trailer only to those stores in higher-priority batches first before the lower-priority ones. In other words, the random set $\J^n \subset \A^n$ consists only of stores from the currently considered high-priority batch. However, because of the min-capacity breach term $-\gamma \sum_{j \in J}b_j$, the function $f(.)$ is not necessarily monotonically increasing, and hence the incremental gain obtained by mapping the next trailer to a store could be negative. In any iteration $n$, if all store-trailer assignments to a high priority batch $\P_p$ yields negative IG, the diminishing property of sub-modular functions eliminates the need for testing trailer assignment to the same batch in subsequent iterations. We proceed to the next lower priority batch $\P_{p-1}$ and choose the store from that batch with the highest IG. This process is repeated until all the eligible stores in the least priority batch $\P_1$ yields negative IG. Once all the batches are processed, we obtain the number of trailers assigned to all the stores. The pseudo-code is given in Sec.~D of the supplementary.
\commentbyKG{
Let $\tilde{\X}$ be the solution obtained from the greedy incremental selection method and $\X^*$ the optimal solution. It is established in ~\cite{Sakaue2020} that for non-negative, non-monotone, submodular functions, the stochastic greedy algorithm guarantees $1/4$ approximation in expectation namely, $\mathbb{E}\left[f(\tilde{\X})\right] \geq \frac{1}{4} (1-\delta)^2 f(\X^*)$ for $0< \delta < 1$, when the condition  $\sum_{j \in \J} R_j \geq 3 \sum_{j \in \J} x_j^*$ is satisfied.
}
\subsection{Transformation to Optimal Transport Problem}
\label{sec:OTtransformation}
The greedy algorithm for the optimal set determination problem measure the computational complexity in terms of the number of calls to a \textit{value oracle} \cite{Buchbinder2017, Kuhnle2019, Sakaue2020}. Given the store-trailer mapping set $\X$, a \textit{value oracle} is a system or an algorithm to compute the value of the submodular function $g(\X)$ in~\eqref{problem:linear_program}. As the greedy algorithm and its variants require repeated evaluation of the incremental gain, \textit{value oracle} is used multiple times in the set selection process. For the solution of our objective in~\eqref{eq:objective} using the stochastic greedy algorithm, the \textit{value oracle} must solve the optimal item allocation problem $\mathbb{P}_1$ in ~\eqref{problem:linear_program}, an instance of LP, once for every candidate store $k$ in the random store set $\J^n$ across multiple iterations $n$. Repeatedly solving multiple instances of LP can take significant amount of time for a large scale supply chain network involving millions of decision variables in $\D$. In this section, we propose to solve the item allocation problem by transforming it to an instance of COT~\cite{Korman2015}, for which the a fast, good quality approximate solution can be obtained via DRM~\cite{Wu2022}.

\begin{theorem}
\label{thm:COT}
Let $P_{ij}^t = \alpha[t] q_{ij}^t$ represent the profit for unit allocation of item $i$ to $j^{th}$ on day $t$. Denote $M_j=Mx_j$ as the maximum possible allocation to store $j$ for the current number of assigned trailers $x_j$. By adding pseudo items and pseudo stores to produce the super-sets $\tilde{\I} \supset I$, and $\tilde{J} \supset J$, and creation of pseudo item inventories $S_i$, max store allocations $M_j$, demands $D_{ij}^t$ and profits $P_{ij}^t$ corresponding to these pseudo items and stores, the LP instance in \eqref{problem:linear_program} can be reduced to an instance of COT:
\begin{equation}
\label{problem:COT}
 \mathbb{P}_2: g(\X) = \max\limits_{\left(d_{ij}^t \right)} \sum_{i,j,t} P_{ij}^t d_{ij}^t
\end{equation}
subject to the constraints:
\begin{alignat*}{2}
    0 \leq d_{ij}^t &\leq D_{ij}^t & \hspace{4mm}&\forall i \in \tilde{\I}, j \in \tilde{\J}, t \in \T_{ij}\\
   \sum_{j \in \tilde{\J}} \sum_{t \in \T_{ij}} d_{ij}^t &= S_i & & \forall i \in \tilde{\I},\\
   \sum_{i \in \tilde{\I}} \sum_{t \in \T_{ij}} d_{ij}^t &= M_j& & \forall j \in \tilde{\J},
\end{alignat*}
satisfying $K = \sum_{i \in \tilde{\I}} S_i = \sum_{j \in \tilde{\J}} M_j$.
\end{theorem}
Sec.~A of the supplementary contains the proof where we also show that the shelf-capacity constraint $s_{ij} \leq C_{ij}$ can be dropped by suitably adjusting the demand values $D_{ij}^t$. The DRM technique is an approximate algorithm and approaches COT by imposing double entropic regularization for both the lower and upper bounds on the transport plan $d_{ij}^t$. The idea is similar to the well known Sinkhorn algorithm~\cite{cuturi13a} popular in OT theory~\cite{peyre19a, villani09a}, where by adding an entropy maximization term $(-\mu \sum_{i,j,t} d_{ij}^t \log\left(d_{ij}^t\right))$, the transport plan can be determined through an  alternate iterative scheme involving simple matrix operations. The computation is efficient as the cubic computational complexity of standard LP solvers are replaced by a linear time solution~\cite{cuturi13a,abid18a}. However because of the additional capacity constraints, in DRM the  matrix-vector multiplication operations of the Sinkhorn algorithm are replaced into finding the unique zero point of several single-variable monotonic functions, determined using any root finding methods~\cite{Press2007}. We refer to~\cite{Wu2022} for detailed description of the algorithm, where the authors argue that the DRM method incurs the same computational complexity as the Sinkhorn iterations. The mathematical formulation is given in Sec.~B of the supplementary for completeness.

We wish to emphasize the following important note. In our framework, the primarily application of DRM is used to quickly identify the candidate store $k_{n+1}$ producing the highest IG given the current store-trailer mapping set $\X^n$, to assign the next trailer $n+1$. After solving the COT instance in~\eqref{problem:COT} for every randomly sampled candidate store $k \in \J^{n+1}$ using DRM, we only require to compare the resultant objective values to determine $k_{n+1}$. As the intermediate item allocation solutions $d_{ij}^t$ are not useful for the entire trailer assignment process described in Sec.~\ref{sec:trailerassignment}, the accuracy required from an approximate method like DRM is enough to \emph{getting the right ordering on the objective values of~\eqref{problem:COT} for different candidate stores.} This allows us to choose a relatively higher value for the parameter $\mu$ controlling the entropy regularization, resulting in faster computation for each COT instance. Such flexibility may not be available when standard LP solvers are used instead.

\section{Experiments}
\label{sec:experiments}
The experiments are run on x$86\_64$ arch., Intel(R) Xeon(R) 2.20GHz CPU, 8 core 64 GB virtual machines, on the \emph{actual data used in the business operations of a global retailer}, for a real supply chain network where a warehouse allocates inventory of $\approx 9000$ items to 158 stores. We use 19 different data sets, each corresponding to a replenishment cycle, to validate our proposed inventory allocation approach. As the data, code and the hyper-parameter settings of $\alpha[t]$, $\beta$ and $\gamma$ are proprietary to the retailer, non-disclosure agreements precludes us from making them publicly available.
\subsection{Set-up}
\label{sec:setup}
As described in Sec.~\ref{sec:INVALS}, the challenge of designing INVALS is in framing the right objective function which makes the inventory allocation plan utility effective and improves business operations. To showcase that our objective in~\eqref{eq:objective} meets these requirements, we perform multiple experiments to evaluate INVALS with respect to the following metrics: (i) Importance of PF for higher labour and trailer utilization, (ii) Need to incorporate store-item priority $q_{ij}^t$ to reduce DOS, (iii) Role played by LTMC trailer cost in improving the trailer utilization and minimizing min-capacity breach, (iv) Experimental validation of Theorem~\ref{thm:COT}, (v) Improvement in computation speed by using DRM method to solve the COT instance in~\eqref{problem:COT} over standard LP solvers, and (vii) Quality of the incremental trailer assignment described in Sec.~\ref{sec:trailerassignment} against the global optimal from off-the-shelf MILP-solver. To this end, we run the following four experiments:
\begin{enumerate}[A.]
    \item Disable PF and no-penalty on LTMC breach ($\gamma=0$).
    \item Enable PF and no-penalty on LTMC breach ($\gamma=0$).
    \item Enable PF with a non-zero penalty for LTMC breach and set $q_{ij}^t = 1$ to enforce constant store-item priority.
    \item Enable all features of INVALS.
\end{enumerate}
 For exp. \textbf{A}, \textbf{B} and \textbf{D}, $q_{ij}^t$ for both coverage and PF periods (if applicable) is set proportional to the demand of item $i$ at store $j$ during time $t$. Higher the value of $q_{ij}^t$, lower the DOS. As the objective function in~\eqref{eq:objective} is maximised, items with lower DOS will be pull-forwarded from future days, and the additional inventory will be held only for few days at the stores leading to smaller holding and maintenance costs. As PF is disallowed in exp. \textbf{A}, for the given data, coverage period need alone is not sufficient to exceed the trailer min-capacity and almost no trailers get dispatched unless $\gamma=0$ to allow this breach. For a high value of $\gamma$ in the no-PF set up, the number of assigned trailers is only $0.64\%$ compared to exp. \textbf{A} and hence we choose to discard this configuration.

For exp. \textbf{D} we ran $4$ different optimization solvers namely, (i) \textbf{Da} (MILP): A direct solution to~\eqref{eq:objective} using off-the-shelf CPLEX~\cite{cplex22} MILP solver to compute global optimum , (ii) \textbf{Db} (LP): Iterative trailer assignment using the GLOP~\cite{glop} solver as the value oracle to solve the item allocation problem in~\eqref{problem:linear_program} and obtain IG for each candidate store, (iii) \textbf{Dc} (COT-LP): Iterative trailer assignment where the GLOP value oracle computes IG by solving the COT instance in~\eqref{problem:COT}, and (iv) \textbf{Dd} (DRM):  Iterative trailer assignment by using the DRM to solve the COT problem and determine IG. Exps. \textbf{A}, \textbf{B} and \textbf{C} are executed with CPLEX to compare at their respective global optimum. For DRM, we set the entropy regularization parameter to $\mu=1$ and used the non-derivative based method such as Brent~\cite{Brent, Dekker} to compute the zero-points of the single variable monotonic functions. Though the Newton's method~\cite{Galantai2000} with its superior quadratic convergence can speed up the computation, it was non-convergent because of small gradient values. The solutions were invariant to small variations around $\mu=1$.


\begin{figure*}[ht]
\captionsetup[subfigure]{}
\centering
\subfloat[]{\includegraphics[width=.48\linewidth]{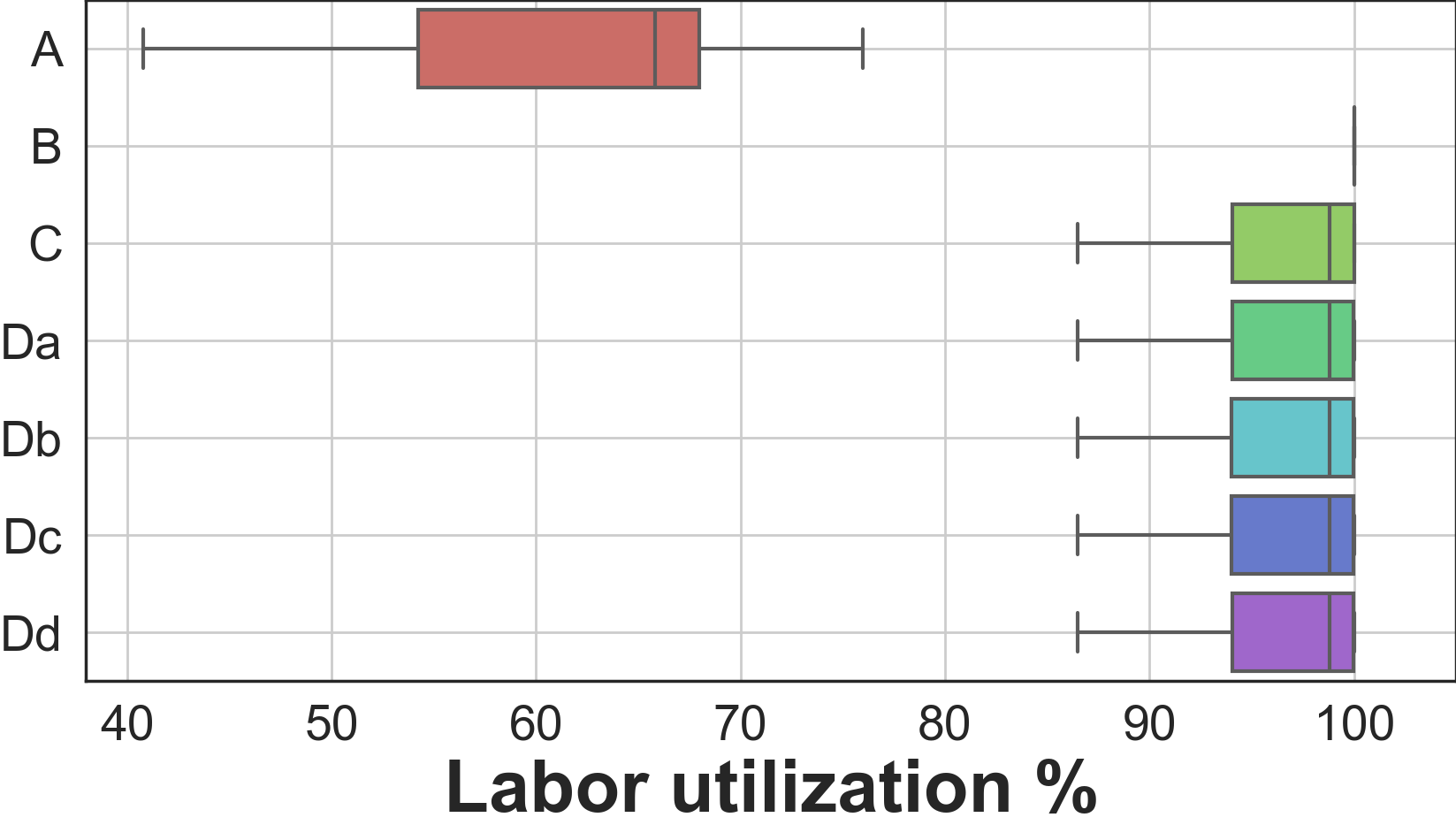}} \hspace{2mm}
\subfloat[]{\includegraphics[width=.48\linewidth]{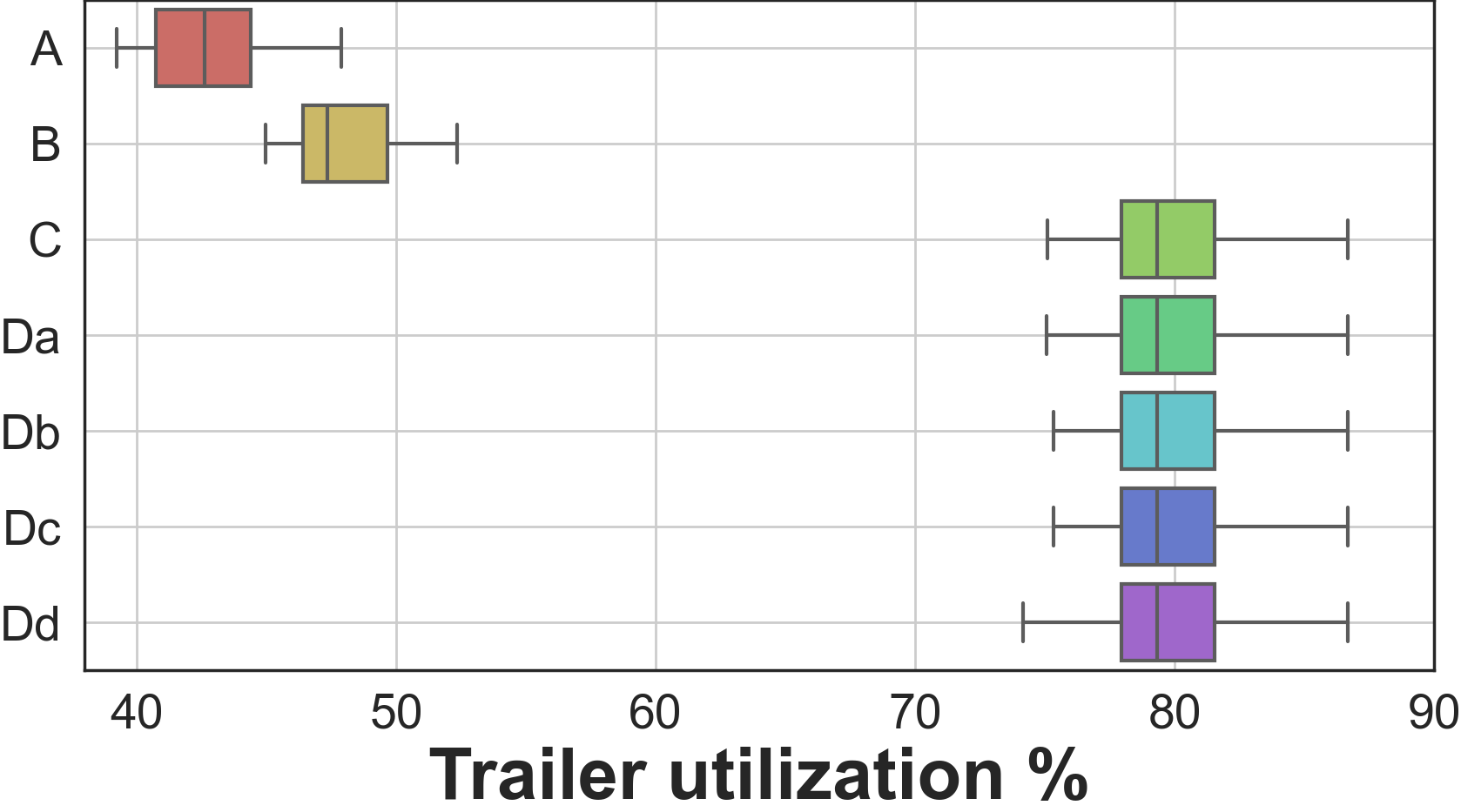}}  \hspace{2mm}
\subfloat[]{\includegraphics[width=.48\linewidth]{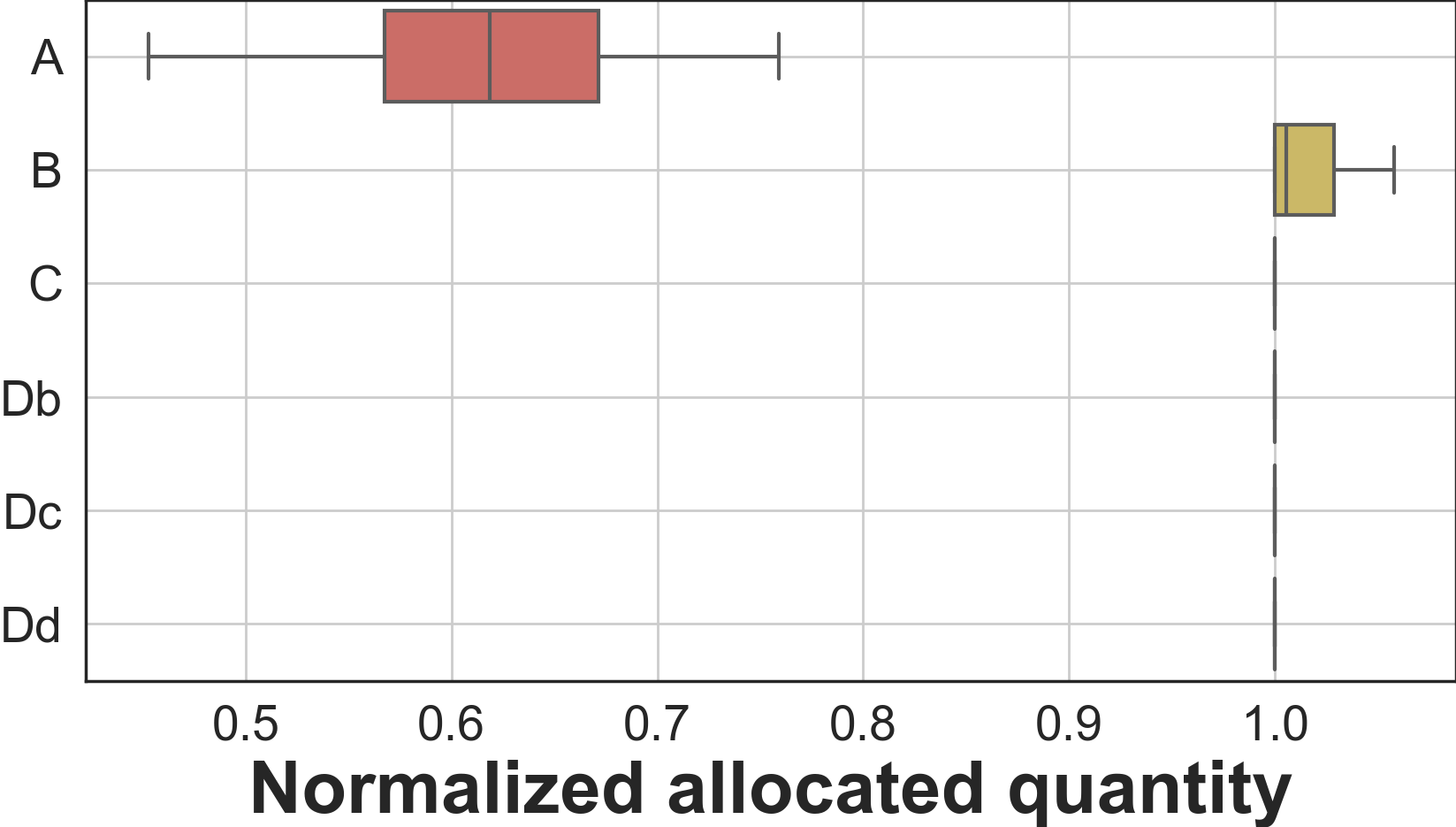}}\hspace{2mm}
\subfloat[]{\includegraphics[width=.48\linewidth]{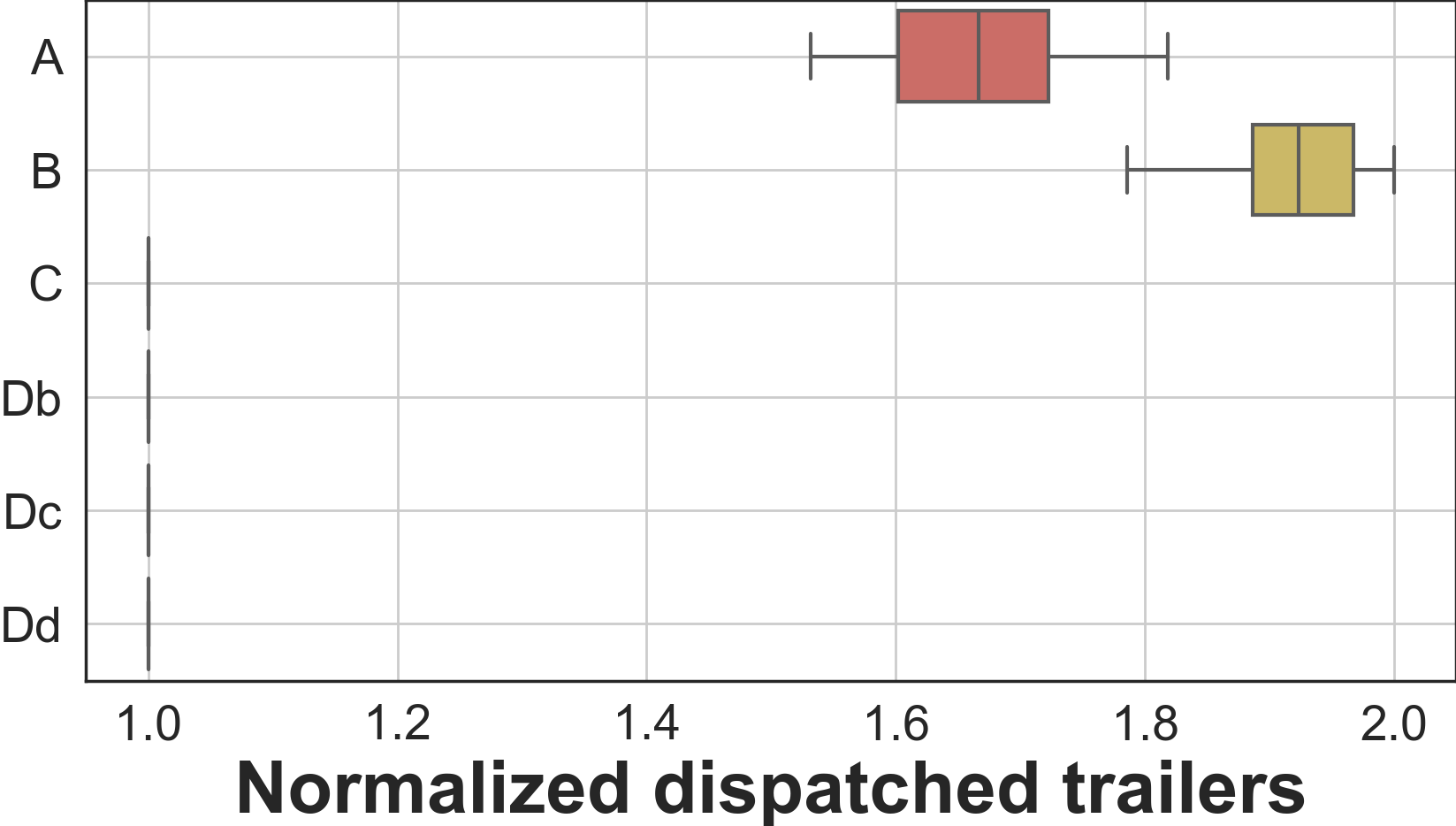}}  \hspace{2mm}
\caption{Disallowing PF in exp. \textbf{A} results in poor trailer and labour utilization, and disabling penalty on LMTC leads to poor trailer utilization in  exp. \textbf{B}. Changing the optimization solver for exp. \textbf{D} between \textbf{Da--Dd} has almost no effect on the trailer and labour utilization, highlighting the quality of our approximate but fast solution approach.}
\label{fig:labtrailer}
\end{figure*}

\begin{figure*}[ht]
\captionsetup[subfigure]{}
\centering
\subfloat[]{\includegraphics[width=.48\linewidth]{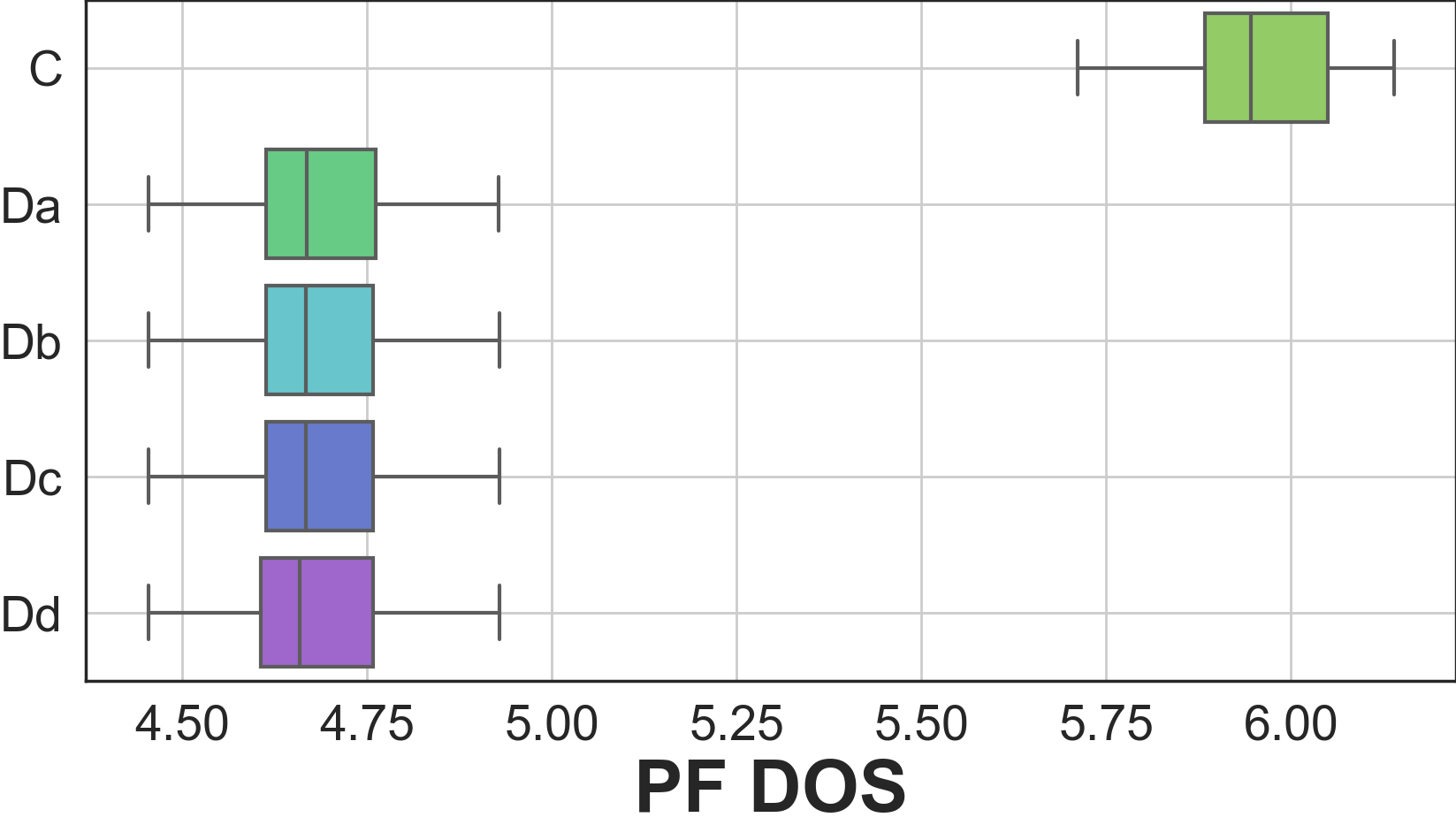}}  \hspace{2mm}
\subfloat[]{\includegraphics[width=.48\linewidth]{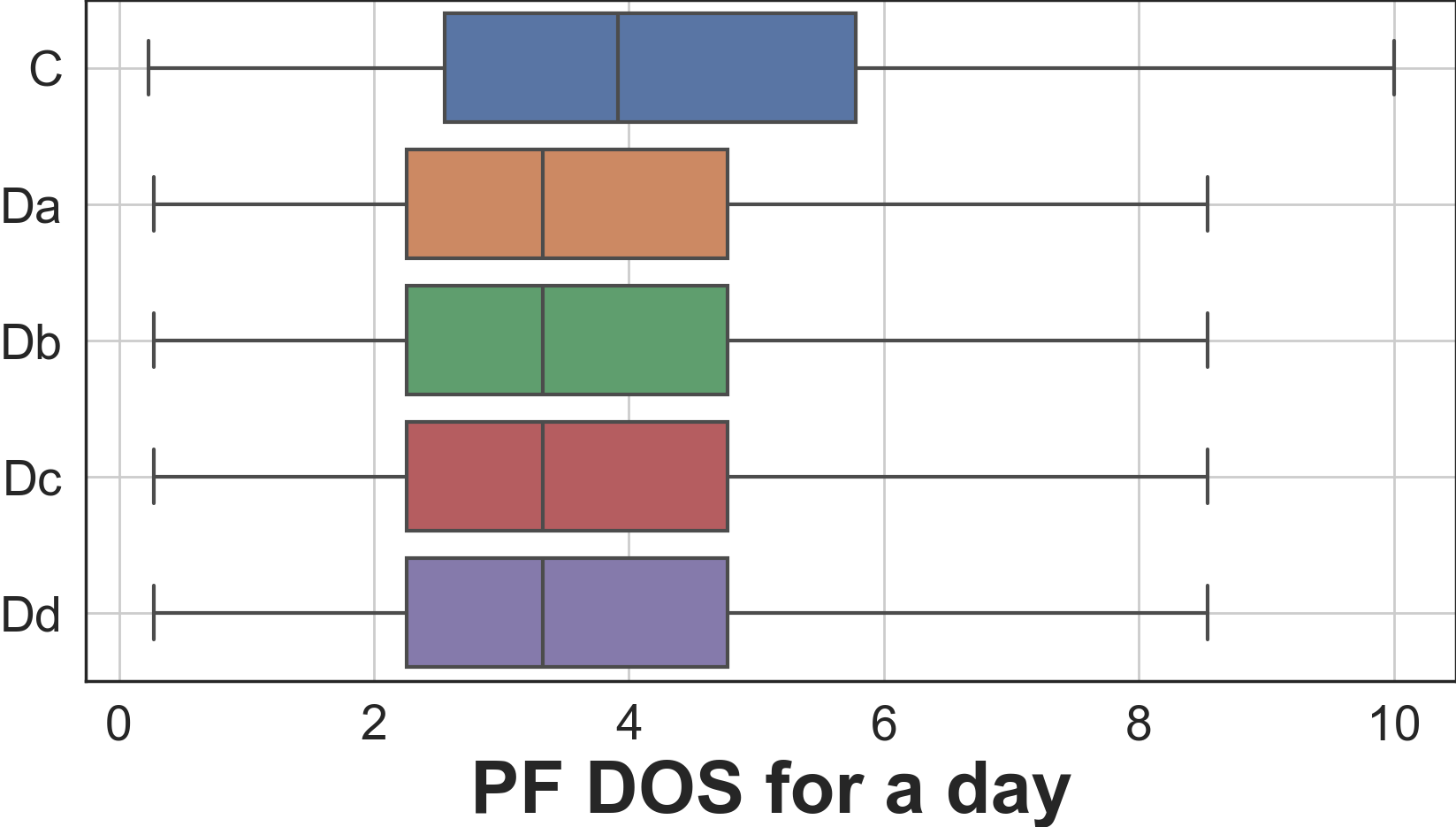}}
\caption{With $q_{ij}^t=1$ in exp. \textbf{C}, the DOS distribution is further to the right of the time axis compared to exps. \textbf{Da--Dd} where it is inversely proportional to store-item DOS. Reduced DOS will decrease the inventory holding costs at the stores and increase savings. Note that the four optimization algorithms in exps. \textbf{Da-Dd} yield almost identical results.}
\label{fig:PFDOS}
\end{figure*}

\begin{figure}[ht]
    \centering
    \includegraphics[width=0.8\linewidth]{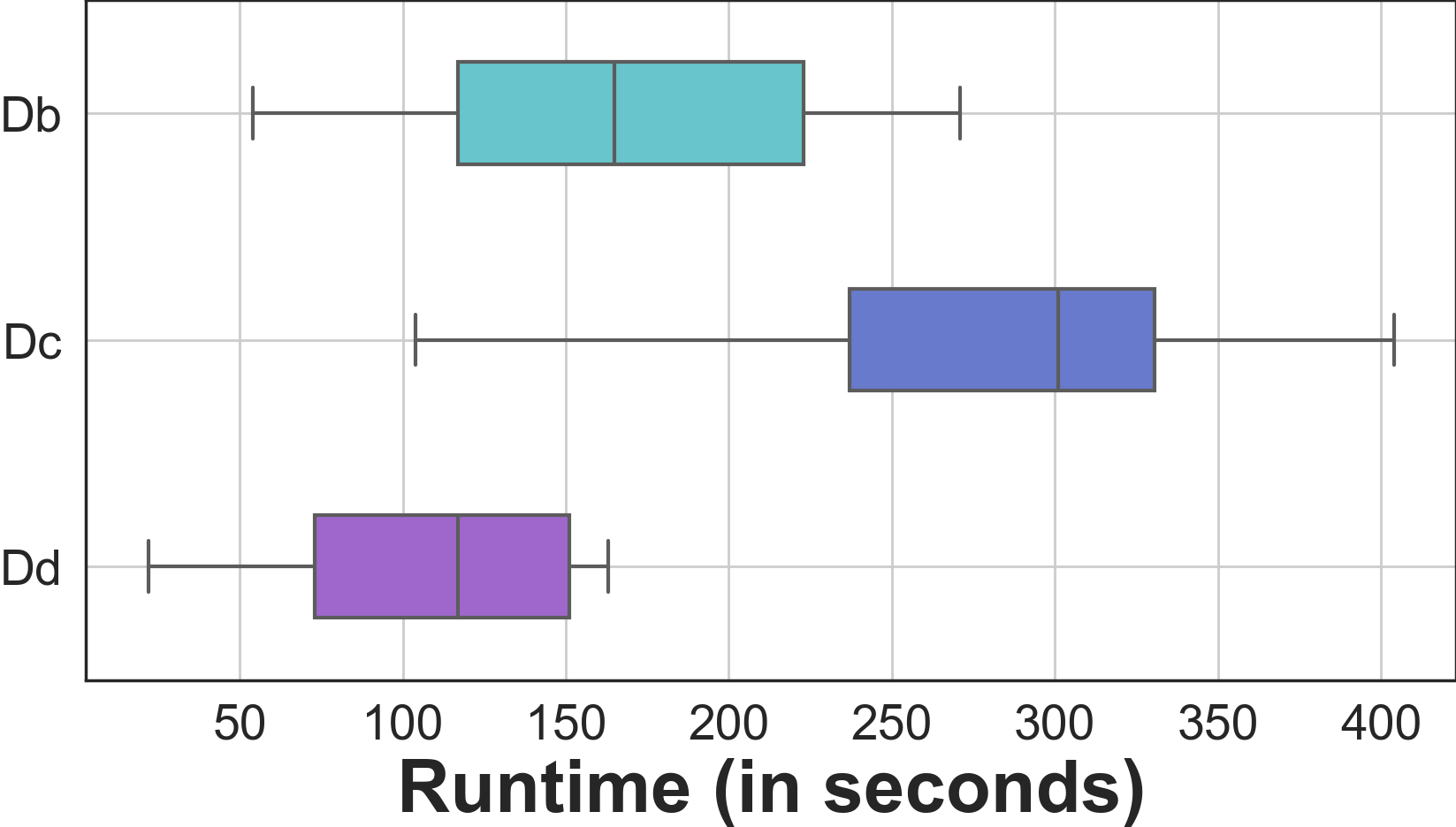}
    \caption{Execution runtime of exps. \textbf{Db}, \textbf{Dc} and \textbf{Dd} for the largest LP instances of problems ~\eqref{problem:linear_program} and \eqref{problem:COT}. DRM in exp. \textbf{Dd} takes $61.16\%$ less time on average than GLOP in exp. \textbf{Dc} for the COT instance in~\eqref{problem:COT}.}
    \label{fig:runtime}
\end{figure}

\subsection{Results}
\label{sec:results}
Figs.~\ref{fig:labtrailer} and \ref{fig:PFDOS} summarize the results of our experiments on multiple data sets from 19 different replenishment cycles. In Figs.~\ref{fig:labtrailer}(a) and \ref{fig:labtrailer}(b) we show the box-plot of the labour and trailer utilization for all the 7 configurations. As expected, exp. \textbf{A}, where the allocation plan is only for the coverage need, makes the least use of the labour and trailer capacities with only $61.73\%$ and $42.61\%$ utilization respectively. Allowing PF in INVALS improves the labour utilization to $96.45\%$. With no-penalty on the min-capacity breach, the trailer utilization in \textbf{B} is only $48.11\%$, and increases substantially to $79.69\%$ in \textbf{C} and \textbf{Da--Dd} where min-capacity is respected because of high $\gamma$ value. In particular, we note that changing the optimization solver for exp. \textbf{D} between \textbf{Da--Dd}, has no detrimental effect on the trailer and labour utilization. This emphasizes the quality of our approximate but fast solution approach. 

In Fig.~\ref{fig:labtrailer}(c) we show the box-plot of the total allocated inventory across all 158 stores. \commentbyKG{normalized with respect to the corresponding allocation from exp. \textbf{Da}.} The values are relative w.r.t. the solution from exp. \textbf{Da} whose allocated inventory is set to $1$ and hence not shown. It is clear that enabling PF increases the size of the inventory plan. Exp. \textbf{B} is expected to have the highest inventory allocation with both PF and min-capacity breach allowed. However with low trailer utilization, exp. \textbf{B} requires $1.92\times$ more trailer than \textbf{Da} as shown in Fig.~\ref{fig:labtrailer}(d) to attain similar size of allocation, underscoring the importance for LTMC trailer cost in~\eqref{eq:objective}. The values in Fig.~\ref{fig:labtrailer}(d) are relative to the total number of trailers allocated in exp. \textbf{Da}. We again see that the trailer assignments from different optimization algorithms are almost identical to each other with the relative values of exps. \textbf{Db--Dd} equalling to $1$.

The importance of store-item priority $q_{ij}^t$ is highlighted in Fig.~\ref{fig:PFDOS}(a) where we show the box-plot of the average DOS computed over the pull-forwarded items across replenishment cycles. The DOS distribution is further to the right of the time axis in exp. \textbf{C} compared to exps. \textbf{Da--Dd}, where $q_{ij}^t=1$ in the former and inversely proportional to store-item DOS in the latter. With an average reduction in DOS by $1.35$ days, the impact on the actual savings in inventory holding costs from INVALS can be expected to be significant. The box-plot of DOS distribution over the pull-forwarded items for one replenishment cycle is shown in Fig.~\ref{fig:PFDOS} (b), where again the DOS values from exp. \textbf{C} are higher compared to exps. \textbf{Da--Dd}.

From Figs. \ref{fig:labtrailer} and \ref{fig:PFDOS} we observe that the four optimization algorithms in exps. \textbf{Da-Dd} yield almost identical results. Even at individual item-store allocation, their values match for all item-stores for $17$ of the $19$ data sets $\approx 90\%$. Recall that once these techniques assign equal number of trailers to each store, they compute the solution for $d_{ij}^t$ from the same (or transformed) LP instance in~\eqref{problem:linear_program}. The solutions from exps. \textbf{Db} and \textbf{Dc} differed from exp. \textbf{Da} for 1 data set due to the stochasticity involved in selecting the set $\J^n$ of candidate stores. We further verified that for each LP instance, invoked for every candidate store $k$ over multiple iterations $n$, exps. \textbf{Db} and \textbf{Dc} gave identical results for all item-store allocations for all data sets, empirically corroborating our Theorem~\ref{thm:COT}. Similarity with the MILP solution validates our approach of assigning and determining store-trailer counts in an incremental fashion. Only for the $2$ data sets, the total number of trailers differ by $2$ between the solutions from exps. \textbf{Da} and \textbf{Dd} . This introduced negligible differences of $0.1\%$ in the trailer utilization and $0.02\%$ in the quantity of inventory allocated between them.

The advantage of employing DRM over LP solver such as GLOP~\cite{glop} is seen in Fig.~\ref{fig:runtime}, where we show the box-plot of the time taken by exps. \textbf{Db}, \textbf{Dc} and \textbf{Dd} for the largest LP instances of problems~\eqref{problem:linear_program} and \eqref{problem:COT}, for each data set. Comparison of absolute times could, in fact, be biased against DRM, as our implementation is pitted against a highly optimized GLOP tool. Nevertheless, we observed that DRM takes $61.16\%$ less time on average than GLOP for the COT instance in~\eqref{problem:COT}. As the addition of pseudo variables increases the instance size in~\eqref{problem:COT} compared to \eqref{problem:linear_program}, the run time for GLOP is longer in \textbf{Dc} than \textbf{Db}. We believe that an optimized implementation and usage of recently proposed root finding algorithms such as~\cite{Kim2017} could further reduce the run time of DRM.

\section{Conclusion}
We studied the problem of inventory allocation from warehouse to stores and proposed a novel forward-looking allocation system, INVALS. We formulated a suitable utility function, which when maximized under different business constraints computes the inventory allocation plan. Exploiting the submodularity of the utility function, we designed an iterative algorithm for finding the number of trailers dispatched to each store from the warehouse, thus eliminating the need to directly solve for integer-valued variables, requiring the usage of MILP solvers. We presented a transformation of the resultant LP problem to a COT instance and leveraged the recently proposed DRM to efficiently compute IG and the selection of the best store to assign the next trailer. We thoroughly investigated the importance of each term in our objective function and analyzed the accuracy of our proposed iterative-approximate algorithm by experimenting with 7 different variants of our framework on 19 data sets. We empirically observed that for 17 $\left(\approx 90\%\right)$ of these data sets, the allocation results from INVALS are identical to the globally optimal MILP solution.

As part of our future work, we would like to develop a similar optimization framework for a multi-echelon system~\cite{Axsater2000}, and determine the end-to-end inventory flow from the supplier through the warehouses to the stores. This could be a challenging initiative, as a large retailer typically purchases products from $100$' suppliers, supplies to $100$' warehouses and generates inventory plans for $1000$' stores. Another avenue of research is to extend our formulation to incorporate store demand variability beyond using only expected demand forecasts. On a general theoretical side, we would like to identify specific sub-structures in instances of LP that make them amenable to be transformed to instances of COT and solved efficiently using methods like DRM.
\appendix
\section{Proof of Theorem 4.1}
Before we proceed with the proof, the following lemma is useful.
\begin{lemma}
\label{lemma:shelfcapacity}
The shelf-capacity constraint $s_{ij} \leq C_{ij}$ for every (item, store) can be dropped altogether by suitable adjustment of the store demands $D_{ij}^t$.
\end{lemma}
\begin{proof}
Let $D_{ij}=\sum\limits_{t \in \T_{ij}}D_{ij}^t$ be to total demand for the $i^{th}$ from store $j$. If $D_{ij} \leq C_{ij}$, then the shelf-capacity constraint is inconsequential as the condition will be implicitly true. Otherwise, consider any two time periods $t_1 < t_2$ and observe that if $d_{ij}^{t_2} >0$, then $d_{ij}^{t_1} = D_{ij}^{t_1}$ for the following reason. Since $\alpha[t]$ is a monotonically decreasing function of $t$, for every (item, store) the objective is maximized by fully allocating the demand for an earlier time period before assigning for later days. This also makes sense from a business perspective, as allocation of any item to a store should first meet the demands of days in the immediate future before utilising for later time periods. Bearing this mind, the demand for the last PF period $D_{ij}^{t_{ij}}$ can be decreased till it reaches zero or $D_{ij}=C_{ij}$, which ever comes earlier. In the former case, we proceed with decreasing the demand for the time period $t_{ij}-1$ and continue the process till $D_{ij} = C_{ij}$.
\end{proof} 

Hence we assume $D_{ij} \leq C_{ij}$ and drop the shelf-capacity constraint for all (item, store). Recall that, the linear program instance to determine the item allocations $\D = \{d_{ij}^t\}$ is invoked for every candidate store $k$ in the randomly chosen set $\J^{n+1}$, for every iteration $n+1$. Let $s_{ij}^n = \sum\limits_t d_{ij}^{t,n}$, $y_j^n$, $b_j^n$ denote the solution of the corresponding decision variables and as before, let $\X^n = \{x_1^n,x_2^n,\ldots,x_J^n\}$ be the number of trailers assigned to each store at the end of iteration $n$.
For the candidate store $k$ in iteration $n+1$ we have $x_k^{n+1} = x_k^n+1$, and for $j \not=k$, $x_j^{n+1} = x_j^n$ and $y_j^{n+1} = y_j^n$. Let $M_j = Mx_j^{n+1}$ $\forall j \in \J$, be the maximum possible allocation to any store using all the trailers assigned to it. 

Without loss of generality, we assume a suitably high value of $\gamma$ where a positive IG for any candidate stores implies that the last trailer assigned to each store is loaded at least to its min-capacity $m$, and all the breach slack variables $b_j$ equal zero at the solution. This is not a restriction to our proof and only simplifies the discussion around the suitable definition of pseudo variables, their corresponding quantities, and the existence of non-empty feasible set. We describe the transformation process below by handling individual constraints in each subsection.

\subsection{Trailer capacity constraints}
\subsubsubsection{Min-capacity constraint}
The negative cost $-\gamma \sum\limits_{j \in \J}b_j$ in the linear objective can be handled by the following additions:
\begin{itemize}
    \item Introduce pseudo item $b$ with inventory $S_b=M$.
    \item Set the profit $P_{bj}^0 = -\gamma$ for the coverage period $t=0$, and $P_{bj}^t = 0$ when $t > 0$ for the PF period, $\forall j \in \J$.
    \item Impose no upper bound on the demand $D_{bj}^0$ for $t=0$, and set $D_{bj}^t = 0$ for $t >0$, $\forall j \in \J$.
\end{itemize}
The straightforward way to handle the constraint $0 \leq b_j \leq m$ would be to set $D_{bj}^0=m$ so that the allocation of the pseudo item $b$ to any store does not exceed the bound $m$. This would entail that $S_b=m*|\J|$. We prove in Section~\ref{sec:feasibility} that for sufficiently high  value of $\gamma$, setting $S_b=M$ and no specific bound on $D_{bj}^t$ is enough to guarantee feasibility and the identification of the best store $k_{n+1}$ having the highest IG$> 0$.

\subsubsection{Max-capacity constraint}
Recall that the max-capacity constraint of $s_j \leq Mx_j^{n+1}$ enforces that the total allocation to the store $j$ across all items does not exceed the available trailer capacity. This condition can be handled as follows:
\begin{itemize}
    \item Introduce pseudo item $z$ with inventory $S_z= (M-m)*|\J|$.
    \item Set profit $P_{zj}^t = 0$, $\forall j \in \J, \forall t \geq 0$. 
    \item Set demand $D_{zj}^0 = M-m$ and $D_{zj}^t = 0$ when $t >0$, $\forall j \in \J$.
    \item Enforce that all trailers are loaded to their maximum capacity $M$ to meet the maximum possible allocation $M_j$.
\end{itemize}
Let $\tilde{\I} = \I \uplus \{b,z\}$. The pseudo inventory $z$ acts as the filler to reach the capacity $M$ in the last trailer assigned to each store. The upper bounds on the demands $D_{zj}^0$ ensure that $z$ is not used more than $M-m$, i.e, the difference between the maximum and minimum trailer capacities.
If item $z$ were to be used, it does not add any profit to the objective function as $P_{zj}^t = 0$. By using the real inventories $i \in \I$, assume that the allocation $s_j^{n+1}$ is deficient in the sense $s_j^{n+1} \leq M\left(x_j^{n+1}-y_j^{n+1}\right) + m$. Then the left over trailer space $M\left(x_j^{n+1}-y_j^{n+1}\right) + m - s_j^{n+1}$ needs to filled with item $b$, incurring a unit cost of $-\gamma$. This is precisely the characteristic of the LP instance in~\eqref{problem:linear_program}. As each $D_{zj}^0 = M-m$, we set $S_z = (M-m)|\J|$.

\subsection{Labour constraint}
For each channel $l \in \L$, we repeat the steps below:
\begin{itemize}
    \item Let $A_l = \sum_{i \in \I_l} S_i$ be the total inventory available across all items belonging to the $l^{th}$ channel. 
    \item Create a pseudo store $h_l$ with need $M_{h_l} = \max(0,A_l- H_l)$, profits $P_{i h_l}^t = 0$, $\forall t \geq 0$, $\forall i \in \tilde{\I}$.
    \item Set demand $D_{i h_l}^t = 0$ $\forall t \geq 0$, $\forall i \in \tilde{\I} \setminus \I_l$.
    \item For $i \in I_l$ impose no upper bound on the demand $D_{i h_l}^0$, and set $D_{i h_l}^t=0$ for $t > 0$.
\end{itemize}
These additions ensure that utmost a sum of $H_l$ quantities of items in $\I_l$ are available for allocation to real stores, to satisfy the warehouse labour constraint $\sum_{i \in \I_l} \sum_{j \in \J} s_{ij} \leq H_l$.

\subsection{Surplus inventory}
Given the fact that we are in iteration $n+1$ where we add a trailer to store $k$, it means the highest IG was positive in iteration $n$. For sufficiently high value of $\gamma$ this positive IG implies that the allocations $s_{ij}^n$ satisfy:
\begin{alignat}{2}
\label{eq:allocation}
\sum\limits_{i \in \I} s_{ij}^n &\geq M(x_j^n-y_j^n)+m, & \hspace{4mm}&\forall j \in \J,\\
\label{eq:satisfylabour}
\sum\limits_{i \in \I_l} \sum\limits_{j \in \J} s_{ij}^n &\leq \min(\sum_{i \in \I_l}S_i, H_l),& & \forall l \in \L \mbox{\hspace{2mm} and}\\
\label{eq:satisfydemand}
d_{ij}^{t,n} &\leq D_{ij}^t & & \forall l \in \I, j \in \J, t \in \T_{ij}. 
\end{alignat}
Recalling that $S_z = (M-m)|\J|$ and $S_b=M$ we have, 
\begin{align*}
\sum\limits_{j \in J}\sum\limits_{i \in \I} s_{ij}^n + S_z + S_b &\geq  \sum\limits_{j \in J} \left[ M(x_j^n-y_j^n)+m+(M-m)\right]\\
&+ M \\
&\geq\left(\sum\limits_{j \not=k} Mx_j^n\right) + Mx_k^{n+1} = \sum\limits_{j \in J} M_j.
\end{align*}
Adding the surplus inventory $\left(\sum\limits_{i \in \I} S_i - \sum\limits_{i \in \I} \sum\limits_{j \in \J} s_{ij}^n\right)$ to both sides we get,
\begin{align*}
    \sum\limits_{i \in \tilde{\I}} S_i &\geq \sum\limits_{j \in J} M_j + \sum\limits_{i \in \I} S_i - \sum\limits_{i \in \I} \sum\limits_{j \in \J} s_{ij}^n\\
    &\geq \sum\limits_{j \in J} M_j + \sum\limits_{i \in \I} S_i - \sum_{l \in \L} \min(\sum_{i \in \I_l}S_i, H_l)\\
    &= \sum\limits_{j \in J} M_j + \sum\limits_{i \in \L} (\sum_{i \in \I_l}S_i) + \sum_{l \in \L} \max(-\sum_{i \in \I_l}S_i, -H_l)\\
    &\stackrel{\rm (a)}{\rm =} \sum\limits_{j \in J} M_j + \sum\limits_{l \in \L}\max\left(0,\sum\limits_{i \in \I_l} S_i - H_l\right) \\
    &= \sum\limits_{j \in J} M_j + \sum\limits_{l \in \L} M_{h_l},
\end{align*}
where the inequality denoted by (a) follows from \eqref{eq:satisfylabour}. The above result establishes that \commentbyKG{when $\gamma$ is high,} the total available inventory is always in surplus to the maximum allocation possible to both real and pseudo stores. To handle this surplus inventory, we introduce one another pseudo store $e$ with need $M_e = \sum_{i \in \tilde{\I}} S_i - \sum_{j \in \J} M_j - \sum_{l \in \L} M_{h_l}$ to absorb real and pseudo item inventories whose allocation is less than their respective $S_i$. For this store $e$:
 \begin{itemize}
     \item Set all profits $P_{ie}^t = 0$, $\forall i \in \tilde{\I}$, $\forall t \geq 0$.
     \item Impose no upper bound on the demand $D_{ie}^0$ for $t=0$, and set $D_{ie}^t =0$ for $t > 0$, $\forall i \in \tilde{\I}$.
 \end{itemize}
 Henceforth, let $\tilde{\J} = \J \uplus \{h_l\}_{l \in \L} \uplus \{e\}$ represent the set of both real and pseudo stores. For small values of $\gamma$, it is possible that the total item inventory is less than then maximum possible store allocation. This setting can be handled in the similar fashion by introducing a pseudo item $e$ instead of a store.

\subsection{Optimal transport problem}
\label{sec:COTproblem}
Armed with this step, the linear program instance for every candidate store $k \in \J^{n+1}$ in the iteration $n+1$ can be re-expressed as:
\begin{equation}
\label{problem:COTappendix}
 \mathbb{P}_2: \max\limits_{\left(d_{ij}^t\right)} \sum\limits_{i,j,t} P_{ij}^t d_{ij}^t
\end{equation}
subject to the constraints:
\begin{alignat}{2}
\label{eq:OTdemandconstraint}
    0 \leq d_{ij}^t &\leq D_{ij}^t & \hspace{4mm}&\forall i \in \tilde{\I}, j \in \tilde{\J}, t \in \T_{ij}\\
\label{eq:OTinventory}
   \sum\limits_{j \in \tilde{\J}} \sum\limits_{t \in \T_{ij}} d_{ij}^t &= S_i & & \forall i \in \tilde{\I},\\
\label{eq:OTneed}
   \sum\limits_{i \in \tilde{\I}} \sum\limits_{t \in \T_{ij}} d_{ij}^t &= M_j& & \forall j \in \tilde{\J}.
\end{alignat}
and satisfying $K = \sum\limits_{i \in \tilde{\I}} S_i = \sum\limits_{j \in \tilde{\J}} M_j$.

\subsection{Feasibility}
\label{sec:feasibility}
We now show that the COT instance always has a non-empty feasible set. \commentbyKG{for sufficiently high $\gamma$ value.}To this end, we prove the following lemma using the inequality~\eqref{eq:allocation} for the allocations $s_{ij}^n$ obtained as the solution to the previous iteration $n$.
 
\begin{lemma}
\label{lem:feasiblesolution}
As before let $s_j^n = \sum_{i \in \I} s_{ij}^n$, $\forall j \in \J$. There exist a solution involving min-capacity breach variables $b_j^{n+1} \geq 0$ satisfying the condition $\sum_{j \in \J} b_j^{n+1}= M$ such that $ s_j^{n} + b_j^{n+1} \geq M\left(x_j^{n+1}-y_j^{n+1}\right)+m, \forall j \in \J$.
\end{lemma}
\begin{proof}
Consider $j \not=k$. Setting $b_j^{n+1}=0$ and using the condition that $s_j^{n} \geq M(x_j^n-y_j^n) +m$ we have $s_j^n+b_j^{n+1} \geq M(x_j^{n+1}-y_j^{n+1}) +m$.
We consider two cases for the store $k$.\\
\noindent \textbf{case 1}: Let $x_k^{n+1} > 1$ implying $y_k^{n+1}=y_k^n$. Setting $b_k^{n+1}=M$ we get $s_k^n+b_k^{n+1} \geq M\left(x_k^n-y_k^n\right)+m+M = M\left(x_k^{n+1}-y_k^{n+1}\right) +m$.\\\\
\noindent \textbf{case 2}: If $x_k^{n+1} = 1$ then $x_k^n=y_k^n=s_k^n=0$, and $y_k^{n+1}=1$. Setting $b_k^{n+1} =M$ we find $s_k^n + b_k^{n+1} \geq M\left(x_k^{n+1}-y_k^{n+1}\right) +m$.\\
The proof follows.
\end{proof}
\noindent We now present a feasible solution for the problem~\eqref{problem:COTappendix} where the variables $b_{j}^{n+1}$ are set to values as stated in Lemma~\ref{lem:feasiblesolution}. 
\begin{itemize}
    \item Set the variables $d_{ij}^{t}$ whose upper bound $D_{ij}^{t} =0$ to $0$.
    \item Set $d_{ij}^{t,n+1} = d_{ij}^{t,n}$, $\forall i \in \I$, $\forall j \in \J$, $\forall t \in \T_{ij}$ satisfying \eqref{eq:satisfydemand} for real items and stores. Then $s_{ij}^{n+1} = s_{ij}^n$, $\forall i \in \I$, $\forall j \in \J$. Letting $d_{zj}^{0,n+1} = M_j - \left(s_j^{n+1} + b_j^{n+1}\right) \leq M-m$ be the allocation of slack item $z$ to store $j$ for day $0$, we have $\sum_{i \in \tilde{\I}} s_{ij}^{n+1} = M_j$, $\forall j \in J$.
    \item From the inequality \eqref{eq:satisfylabour} it follows that the allocations $s_{ij}^{n+1}$ also respect the labour constraint. Then for each channel $l \in \L$ with the corresponding store need $M_{h_l} > 0$, there exists non-negative variables $d_{ih_l}^{0,n+1} \leq S_i - \sum_{j \in \J} s_{ij}^{n+1}$ for $i \in \I_l$, such that  $\sum_{i \in \I_l}s_{ih_l}^{n+1} = M_{h_l}$. For those categories whose need $M_{h_l}=0$, set $d_{ih_l}^{0,n+1}=0$ $\forall i \in \I_l$.
    \item Finally, set $d_{ie}^{0,n+1} = S_i - \sum_{j \in \J} s_{ij}^{n+1}-\sum_{l \in \L}s_{ih_l}^{n+1}$, $\forall i \in \tilde{\I}$ assigning the left over inventories to pseudo store $e$. By construction of the store need $M_e$ we have $\sum_{i \in \tilde{\I}} s_{ie}^{n+1}=M_e$. 
\end{itemize}  
For the best store $k_{n+1}$ having the highest IG$>0$, the solution will satisfy the condition $s_j^{n+1} >= M(x_j^{n+1}-y_j^{n+1})+m$, for all real stores $j \in \J$. Item $z$ will used as a filler in the last trailer so that the store need of $M_j = Mx_j^{n+1}$ is met. For all $j \in \J$ allocation $s_{bj}^{n+1} = 0$, and for the dummy store $e$, $s_{be}^{n+1} = M$ to absorb the inventory of the pseudo item $b$. As $P_{be}^t = 0$, no negative cost is added to the objective~\eqref{problem:COTappendix} from item $b$.

\section{Mathematical Formulation of the DRM}
For the COT instance in Section~\ref{sec:COTproblem} we can replace $S_i \leftarrow S_i / K$, $M_j \leftarrow M_j / K$, and $D_{ij}^t \leftarrow D_{ij}^t / K$ to make the total inventory flow in the network to sum up to 1. The resultant solution $d_{ij}^t$ can be scaled back by $K$ to obtain the true inventory allocation. The COT instance differs from the classical OT because of the presence of a point-wise upper bound $d_{ij}^t \leq D_{ij}^t$ on each element of the transport plan, limiting the mass transported between each pair of item (source) and store (sink). Further, the presence of multiple days $t$ adds an extra third time dimension to the COT problem. However, as no constraints are present along the time axis because the shelf-capacity constraints can be omitted as shown in Lemma~\ref{lemma:shelfcapacity}, the complexity of problem~\eqref{problem:COTappendix} is the same as the setting studied in~\cite{Wu2022}.

Define the sets 
\begin{align*}
    \M &= \{(i,j,t) : D_{ij}^t > 0 \text{ or unspecified} \}, \\
    \N &= \{(i,j,t): (i, j, t) \in \M \text{ \& } D_{ij}^t \text{ is finite}\}
\end{align*}
For $(i,j,t) \notin \M$, $d_{ij}^t=0$ as the demand $D_{ij}^t=0$. By introducing entropic regularization for both the lower and upper demand bounds on the allocation variables $d_{ij}^t$ and the addition of Lagrange parameters, the DRM method in~\cite{Wu2022} considers the modified objective:
\begin{align}
\label{eq:entreg}
     \O&= \sum\limits_{\M} P_{ij}^t d_{ij}^t - \mu\sum\limits_{\M}d_{ij}^t \ln\left(d_{ij}^t\right) \nonumber \\
     &-\mu\sum\limits_{\N} \left(D_{ij}^t-d_{ij}^t\right) \ln\left(D_{ij}^t-d_{ij}^t\right) \nonumber \\
     &- \sum_{i}\alpha_{i}\left(\sum_{jt} d_{ij}^t - S_i\right)- \sum_{j}\beta_j \left(\sum_{it} d_{ij}^t - M_j\right).
\end{align}
By taking partial derivatives with respect to $d_{ij}^t$ and setting $\frac{\partial \O}{\partial d_{ij}^t} = 0$, we get
\begin{equation*} 
\label{cot:gamma}
    d_{ij}^t =
    \begin{cases}
        D_{ij}^t \left[1- \frac{1}{1 + \phi_{i}K_{ij}^t\psi_{j}}\right],& \text{if } (i,j,t)\in \N\\
        \frac{\phi_{i}K_{ij}^t\psi_{j}}{\mathrm{e}},& \text{if } (i,j,t)\in \M-\N
    \end{cases}
\end{equation*}
where $\phi_i = \mathrm{e}^{-\frac{\alpha_i}{\mu}}$, $\psi_j = \mathrm{e}^{-\frac{\beta_j}{\mu}}$ and $K_{ij}^t = \mathrm{e}^{\frac{P_{ij}^t}{\mu}}$.
The problem reduces to identifying $|\tilde{I}|$ values $\phi_i$ and $|\tilde{J}|$ values $\psi_j$, for a total of $|\tilde{I}|+|\tilde{J}|$ constants, such that the supply and demand constrains in \eqref{eq:OTinventory} and \eqref{eq:OTneed} are simultaneously meet. DRM determines these constants via an alternative iteration scheme similar to the Sinkhorn algorithm, as described below. 

For some large constant $\zeta$, we can set $D_{ij}^t = \zeta$ for $(i,j,t) \in \M\text{ - }\N$ and without loss of generality assume $\N = \M$. Setting $\frac{\partial \O}{\partial \alpha_i} = \frac{\partial \O}{\partial \beta_j} = 0$ of the objective in~\eqref{eq:entreg}, we have
\begin{align*}
 S_i &= \sum\limits_{jt}\left(D_{ij}^t - \frac{D_{ij}^t}{1 + \phi_{i}K_{ij}^t\psi_{j}} \right)\\
  M_j &= \sum_{it}\left(D_{ij}^t - \frac{D_{ij}^t}{1 + \phi_{i}K_{ij}^t\psi_{j}}\right).
\end{align*}
Consider an alternative minimization scheme where $\{\psi^m\}$ are the values at the end of iteration $m$. Then each $\phi_i^{m+1}$ can be determined as the  \emph{unique zero point} of its corresponding single-variable monotonic function $g_i(\phi_i)$ given by:
\begin{equation*}
    \label{eq:gphii}
    g_i(\phi_i) = S_i - \sum\limits_{jt} D_{ij}^t + \sum\limits_{jt} \frac{D_{ij}^t}{1 + \phi_{i}K_{ij}^t\psi_{j}^{m}}.
\end{equation*}
Observe that $g_i(0) > 0$ and $g_i(\infty) < 0$ and $\phi_i^{m+1}$ is the point where $g_i\left(\phi_i^{m+1}\right) = 0$. Similarly, given $\{\phi_i^m\}$ each $\psi_j^{m+1}$ is obtained as the unique zero point of its corresponding single-variable monotonic function $h_j(\psi_j)$ given by:
\begin{equation*}
    \label{eq:hphij}
    h_j(\psi_j) = M_j - \sum\limits_{it} D_{ij}^t + \sum\limits_{it} \frac{D_{ij}^t}{1 + \phi_{i}^m K_{ij}^t\psi_{j}}.
\end{equation*}
The zero points can be computed using any root finding algorithm in one-variable such as Newton's, Bisection, Interpolate Truncate and Project, Dekker's, Brent's, etc.~\cite{Press2007}. Such an alternating scheme is proven to converge in~\cite{Wu2022}.

\section{Submodularity of the Item Allocation Linear Program}
\begin{definition}[Incremental Gain]
For a set function f($\cdot$), a subset $\mathcal{A}\subseteq \mathcal{V}$ and an element $i\in \mathcal{V}$, the incremental gain is defined as:
\begin{equation*}
\label{eq:incrementalgain}
f_\mathcal{A}\left(i\right)=f\left(\mathcal{A}\cup\{i\}\right)\ -f\left(\mathcal{A}\right)	    
\end{equation*}
\end{definition}

\begin{definition}[Modularity, Submodularity, Non-negativity]
Consider any two sets $\mathcal{A}\subseteq\mathcal{B}\subseteq\mathcal{V}$. A set function $f(\cdot)$ is submodular iff for any $i\notin\mathcal{B}$, $f_\mathcal{A}(i)\geq f_\mathcal{B}(i)$, non-negative when $f(\mathcal{A}) \geq 0$ and modular when $f_{\A}(i)$ is independent of $\A$.
\end{definition}
Recall our representation $\X = \{x_1,x_2,\ldots,x_{|\J|}\}$ and the slight abuse of notation $\Z = \X \uplus \{k\}$ as the set of store-trailer mappings with $z_k = x_k+1$, and $z_j = x_k, k\not=j$. As $h_{\X}(k) = p_k$, it is a modular function.

We now study the submodularity of the item allocation linear program $g(\X)$ defined in~\eqref{problem:linear_program}. Consider two sets $\X$ and $\Y = \{y_1,y_2,\ldots,y_{|\J|}\}$ where $\X \subseteq \Y$ in the sense that $x_j \leq y_j, \forall j$. To prove the submodular property, we require to show that the diminishing returns property 
\begin{equation}
\label{eq:diminishingreturns}
g(\X \uplus \{k\})-g(\X) \geq g(\Y \uplus \{k\})-g(\Y)
\end{equation}
holds after addition of a new trailer to any store $k$, given the current store-trailer assignments. Under a much restricted setting of $R_j=1$, $|\T_{ij}|=1$, uniform item inventory of $S_i = \kappa, \forall i$ and no constraints on labour $H_l$, trailer min-capacity $m$, and in particular no point-wise capacity constraint $D_{ij}^t$ limiting the allocation between each pair of item and store, the LP instance in~\eqref{problem:linear_program} reduces to a \emph{relaxed assignment problem} which is proven to be submodular in~\cite{Kulik2019}. 

For our generalized formulation we do not yet have a formal mathematical proof of this property. It appears to be an arduous task primarily due to the individual capacity constraints $d_{ij}^t \leq D_{ij}^t$. To see this, consider the interesting case where the incremental gain $g_{\Y}(k) = g(\Y \uplus \{k\})-g(\Y)$ is positive. As discussed before, a positive IG for our parameter setting of high $\gamma$ would imply that at the solution $\left(\D^*_{\Y \uplus \{k\}}, \B^*_{\Y \uplus \{k\}}\right)$ corresponding to the set $\Y \uplus \{k\}$ for the LP instance in~\eqref{problem:linear_program}, all trailers are loaded to their min-capacities, and the allocations $s_{ij}^{\Y \uplus \{k\}}$ satisfy~\eqref{eq:allocation}. The min-capacity slack variables $b_j^{\Y \uplus \{k\}} =0, \forall j \in \J$. If the inventory allocation plan for a super-set $\Y$ does not involve min-capacity breaches, then the allocation plan for any subset $\X \subseteq \Y$ should also respect the min-capacity constraint with $b_j^{\X}=0, \forall j$. Consider the equivalent COT instance stated in Sec.~\ref{sec:COTproblem}. Then the solution for each of the four store-trailer assignment sets $\X$, $\X \uplus \{k\}$, $\Y$, and $\Y \uplus \{k\}$ will have the structure that $s_{be} = M$, and $s_{bj}=0, \forall j \in \J$, with the inventory of pseudo item $b$ allocated to the pseudo store $e$ in entirety as detailed in Sec.~\ref{sec:feasibility}. It becomes redundant to add the pseudo inventory $b$ and its corresponding pseudo quantities in the transformation to the COT instance. Then the profit matrix with entries $P_{ij}^t$ in ~\eqref{problem:COTappendix} will all be non-negative. Further, the total item inventory of $K=\sum_{i \in \tilde{\I}} S_i$ is a constant, independent of the store-trailer mapping set $\X$. It is only the maximum store allocation values $M_j$ and the pseudo store need $M_e$ that vary with the COT instances. Our problem structure gradually begins to take the shape of the assignment problem studied in~\cite{Kulik2019}. It is the presence of the capacity constraints $d_{ij}^t \leq D_{ij}^t$ that makes the COT instance fundamentally different, engendering complexity to the proof of submodularity.

\begin{figure}
    \centering
    \includegraphics[width=8.5cm]{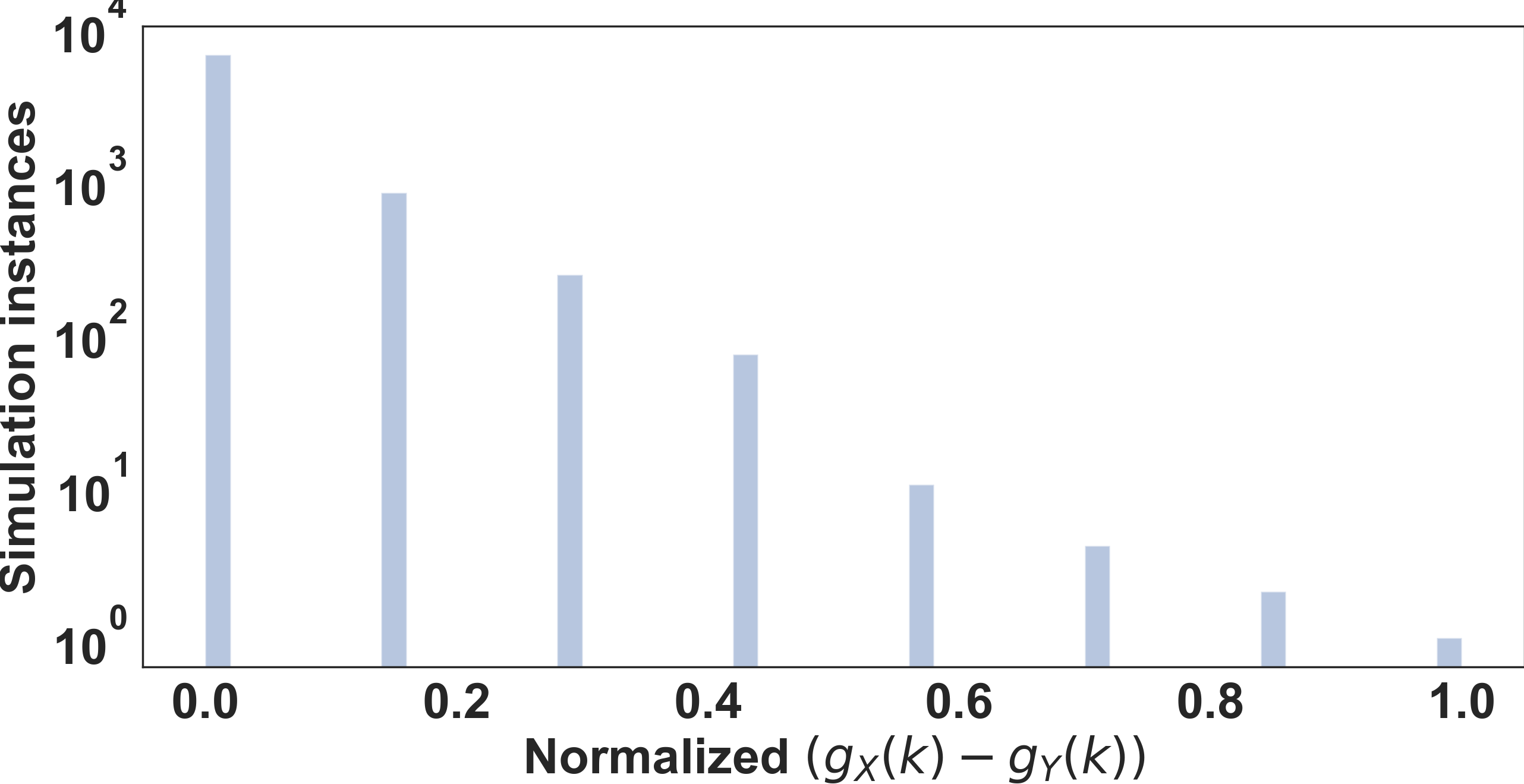}
    \caption{Distribution of normalised $g_{\X}(k) - g_{\Y}(k)$}
    \label{fig:distribution-of-positive-deltas}
\end{figure}

However, we firmly believe that $g(\X)$ is submodular based on our experimental validations. To this end, we created $760$ sets of store-trailer mappings $\X \subset \Y$ such that the LP instance in~\eqref{problem:linear_program} always had a feasible solution for $\Y$ (and hence for $\X$). For each pair of $(\X, \Y)$, we considered $10$ different choices of candidate stores $k$ to assign the next trailer while ensuring that $\Y \uplus \{k\}$ had a feasible solution. For each triplet $(\X, \Y, k)$, we computed the incremental gains $g_{\X}(k)$ and $g_{\Y}(k)$ and generated a distribution plot of the difference $g_{\X}(k)-g_{\Y}(k)$ consisting of $7600$ data-points. The normalised plot is shown in Fig.~\ref{fig:distribution-of-positive-deltas} where the differences are made relative to the largest value of $g_{\X}(k)-g_{\Y}(k)$ and the latter set to 1. \emph{We observe that every difference is non-negative}; thus empirically corroborating the diminishing returns property~\eqref{eq:diminishingreturns} of the function $g(.)$. 
\commentbyKG{
\begin{figure}
    \centering
    \includegraphics[width=8.5cm]{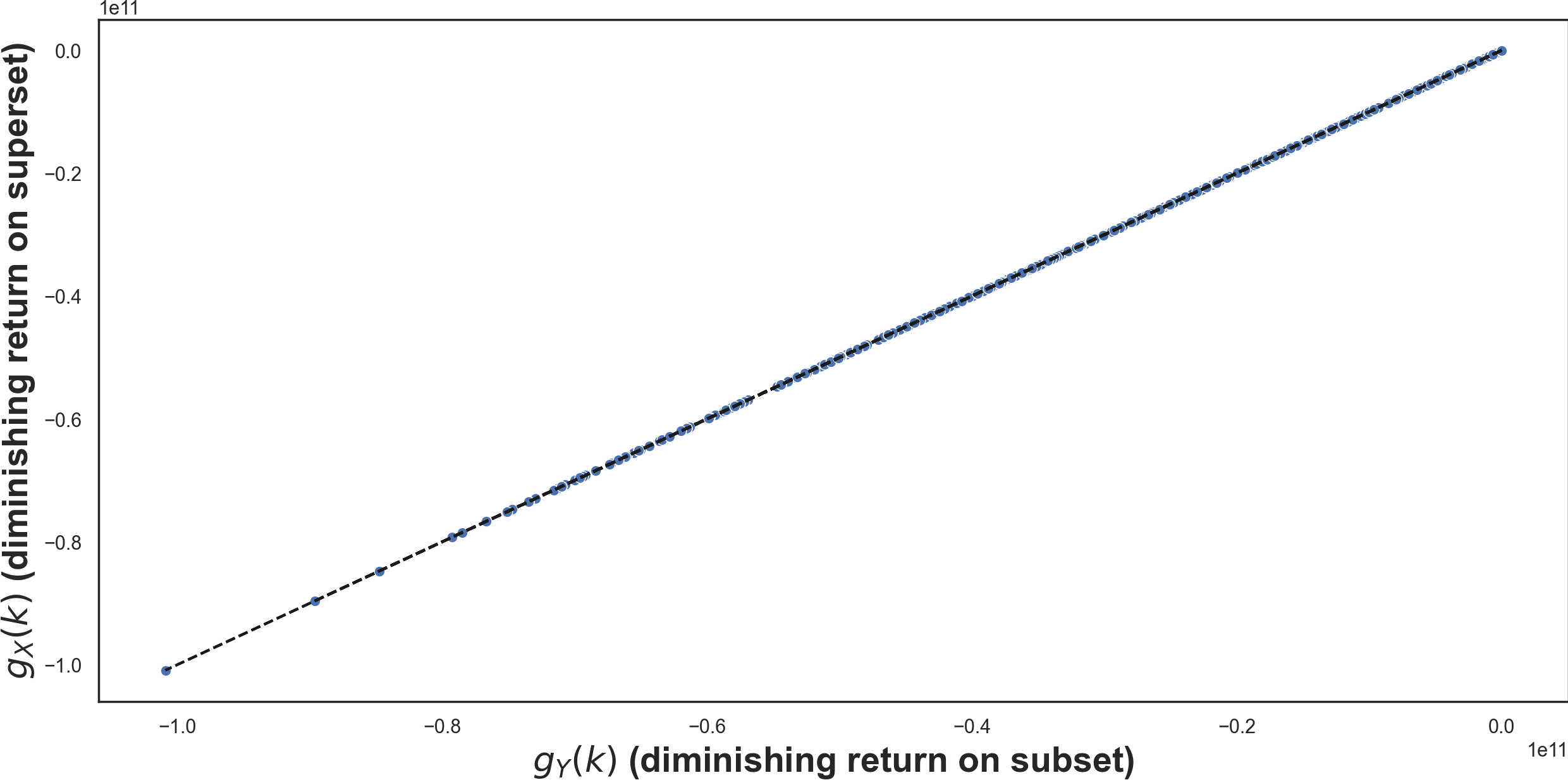}
    \caption{Submodular property}
    \label{fig:diminishingreturns}
\end{figure}
}

\begin{algorithm*}[ht!]
\caption{Batched stochastic greedy algorithm}
\label{algo:stochastic_reedy}
\begin{algorithmic}[1]
\State {Set $x_j=0, \forall j$}
\State {Choose $\rho \in (0,1]$} \Comment{Store selection probability}
\For {$p$ in decreasing store priority}
    \State $\X = \{x_1,x_2,\ldots,x_{|\J|}\}$
    \State $\P_p \gets \{j \in \J: p_j=p\}$ 
    \While {$\P_p \not= \emptyset$}
        \State $(k, \mathcal{F})$ = NextBestStore$(\X, \P_p, \rho)$ \Comment{Determine the best store to assign next trailer}
        \State $\P_p\gets \P_p \setminus \mathcal{F}$ \Comment{Remove stores with negative or zero IG}
        \If {$k \not= \varnothing$} 
            \State $x_{k} \gets x_{k}+1$
            \If {$x_k=R_k$} \Comment{Planned trailer constraint met}
                 \State $\P_p \gets \P_p \setminus \{k\}$
            \EndIf
        \EndIf
    \EndWhile
\EndFor
\State \Return $\X$
\State {\textbf{Output}: $\X$-Final set of store-trailer assignments}
\end{algorithmic}
\end{algorithm*}

\begin{algorithm*}[ht!]
\caption{NextBestStore $(\X, \P, \rho)$}
\label{algo:next_store}
\begin{algorithmic}[1]
\State {\textbf{Input}: $\X$-Current store-trailer assignments, $\P$-Candidate stores, $\rho$-Store selection probability}
\State {Set $\mathcal{F} \gets \emptyset$, $k_n \gets \varnothing$}
\State $\J^n \gets$ select stores with probability $\rho$ from $\P$ \Comment{Random selection of candidate stores}
\For{$k \in \J^n$}
    \State  Calculate $g_{\X}(k) \gets g\left(\X \uplus \{k\}\right) - g\left(\X\right)$ \Comment{Increment gain}
     \If{$g_{\X}(k) \leq 0$}
        \State $\mathcal{F} \gets \mathcal{F} \uplus \{k\}$ \Comment{List of stores with negative or zero IG}
    \EndIf
\EndFor
\State $k_n \gets \argmax_{k \in \J^n} g_{\X}(k)$ \Comment{Store with highest IG}
\If{$g_{\X}(k_n) \leq 0$}
    \State $k_n \gets \varnothing$
\EndIf
\State \Return $(k_n, \mathcal{F})$
\State {\textbf{Output}: $k_n$-Next best store, $\mathcal{F}$-Set of stores with negative or zero IG}
\end{algorithmic}
\end{algorithm*}
\section{Pseudo-code of the Incremental Trailer Assignment Method}
Below we present the pseudo-code of the batched iterative algorithm to determine the best store for assigning the next trailer. To keep it simple, we have not incorporated the \emph{laziness} advantage of greedy algorithms~\cite{minoux78a, leskovec07a} in our description. It can be included using a max-heap data structure and further constraint the stores randomly selected in step 3 of Algorithm~\ref{algo:next_store}.

\bibliography{INVALS}
\end{document}